\documentstyle[amssymb,epsfig,11pt]{article}
\newtheorem{Lemma}{Lemma}[section]
\newtheorem{Theorem}[Lemma]{Theorem}
\newtheorem{Proposition}[Lemma]{Proposition}
\newtheorem{Corollary}[Lemma]{Corollary}
\newtheorem{Remarks}[Lemma]{Remarks}
\newtheorem{Definition}[Lemma]{Definition}
\newtheorem{Assumption}[Lemma]{Assumption}
\newtheorem{Example}[Lemma]{Example}
\newtheorem{Notation}[Lemma]{Notation}
\def\theequation{\arabic{section}.\arabic{equation}}
\newcommand{\isep}{\itemsep -0.6mm}

\newcommand{\rr}    {\rightarrow}

\newcommand{\ts}{{\tilde{s}}}
\newcommand{\tf}{{\tilde{f}}}
\newcommand{\tx}{{\tilde{x}}}
\newcommand{\ty}{{\tilde{y}}}

\newcommand{\us}     {\underline{s}}
\newcommand{\os}     {\overline{s}}
\newcommand{\UP}     {\underline{P}}
\newcommand{\OP}     {\overline{P}}

\newcommand{\HD}    {\mbox{${\rm dim}_H$}}
\newcommand{\dimH}    {{\rm dim}_H}
\newcommand{\boxinf}    {{\rm dim}_B}
\newcommand{\boxsup}    {{\rm dim}^B}

\newcommand{\re}    {{\rm re}}
\renewcommand{\Re}    {{\rm Re}}

\newcommand{\Aut}    {{\rm Aut}}
\newcommand{\Area}    {{\rm Area}}

\newcommand{\crit}    {{\rm crit}}

\newcommand{\diam}    {{\rm diam}}
\newcommand{\bfh}   {{\bf h}}
\newcommand{\bfone}   {{\bf 1}}
\newcommand{\Int}    {{\rm Int}}

\newcommand{\half} {\frac{1}{2}}

\newcommand{\chiLz}{\chi_{\mbox{}_{\Lambda_0}}}
\newcommand{\chiU}{\chi_{\mbox{}_{U}}}
\newcommand{\chiUL}{\chi_{\mbox{}_{U \cap \Lambda_0}}}

\newcommand{\calA}{{\cal A}}
\newcommand{\CCC}{{\cal C}}
\newcommand{\dC}{{{d_{\cal C}}}}
\newcommand{\D}{{\cal D}}
\newcommand{\EE}{{\cal E}}
\newcommand{\dE}{{d_{\cal E}}}
\newcommand{\FF}{{\cal F}}
\newcommand{\tFF}{\tilde{{\cal F}}}

\newcommand{\LL}{{\cal L}}
\newcommand{\calM}{{\cal M}}
\newcommand{\OO}{{\cal O}}
\newcommand{\PP}{{\cal P}}
\newcommand{\calP}{{\cal P}}

\newcommand{\tbff}{{\tilde{{\bf f}}}}
\newcommand{\bff}{\mbox{$\bf f$}}
\newcommand{\bomega}{\mbox{$\bf \omega$}}

\newcommand{\CC}{\mbox{$\Bbb C$}}
\newcommand{\CCs}{{\Bbb C}}
\newcommand{\HCC}{{\widehat{\Bbb C}}}
\newcommand{\tinyD}{{\mbox{}_{\Bbb D}}}
\newcommand{\htinyD}{{\mbox{}_{\widehat{\Bbb D}}}}
\newcommand{\DDs}{{\Bbb D}}
\newcommand{\DD}{{\Bbb D}}
\newcommand{\EEE}{\mbox{$\Bbb E$}}

\newcommand{\RR}{\mbox{$\Bbb R$}}
\newcommand{\RRs}{{\Bbb R}}
\newcommand{\NN}{\mbox{$\Bbb N$}}
\newcommand{\NNs}{{\Bbb N}}

\newcommand{\AutDD}    {{\rm Aut(\DD;\diag \DD)}}

\newcommand{\bara}{\overline{a}}
\newcommand{\barb}{\overline{b}}
\newcommand{\barf}{\overline{f}}

\newcommand{\bart}{{\overline{t}}}
\newcommand{\baru}{\overline{u}}
\newcommand{\barx}{\overline{x}}
\newcommand{\bary}{\overline{y}}
\newcommand{\barz}{\overline{z}}
\newcommand{\barw}{\overline{w}}
\newcommand{\bardz}{d\overline{z}}

\newcommand{\barDD}{{\overline{\DD}}}

\newcommand{\barP}{\overline{P}}
\newcommand{\barR}{\overline{R}}
\newcommand{\barU}{\overline{U}}

\newcommand{\UbarU}{U\times \overline{U}}

\newcommand{\barphi}{\overline{\phi}}
\newcommand{\barpsi}{\overline{\psi}}

\newcommand{\hatd}{\widehat{d}}
\newcommand{\hatg}{g_\htinyD}

\newcommand{\hatpsi}{\widehat{\psi}}
\newcommand{\hatr}{\widehat{r}}
\newcommand{\hatx}{\widehat{x}}
\newcommand{\haty}{\widehat{y}}

\newcommand{\hatD}{\widehat{\DD}}

\newcommand{\hatR}{\widehat{R}}
\newcommand{\hatU}{\widehat{U}}

\newcommand{\HK}{{\widehat{K}}}
\newcommand{\KD}{{\widehat{K}}_\Delta}

\newcommand{\hatf}{{\widehat{f}}}

\newcommand{\Cl}{\mbox{Cl}}
\newcommand{\len}{\mbox{len\,}}
\newcommand{\diag}{\mbox{diag\ }}
\newcommand{\diagK}{\mbox{diag}(K)}

\newcommand{\Halmos}    {\ \raisebox{0.3ex}  {\framebox[0.9ex]{
			   \rule[0ex]{0ex}{0.5ex}
			    }}}

\begin{document}
\title{ On the dimensions of conformal repellers.
Randomness and parameter dependency}
\author{Hans Henrik Rugh\\
    University of Cergy-Pontoise, France}
\date {\today}
 \maketitle
\begin{abstract}
Bowen's formula relates the Hausdorff dimension of a conformal repeller
to the zero of a `pressure' function. 
We present an elementary, self-contained
proof which bypasses measure theory and
the Thermodynamic Formalism to show that  
Bowen's formula holds for $C^1$
conformal repellers. 
We consider time-dependent conformal repellers
 obtained as invariant
subsets for 
sequences of conformally expanding maps within a suitable class.
We show that Bowen's formula generalizes to such a  repeller 
and that if the sequence is picked at random
then the Hausdorff dimension of the repeller
almost surely agrees with its upper and lower
Box dimensions and is given by a natural
generalization of Bowen's formula.
For a random uniformly hyperbolic Julia set 
on the Riemann sphere we show that if the family of maps
and the probability law depend real-analytically
on parameters then so does its almost sure Hausdorff dimension.
\end{abstract}

\section{Random Julia sets and their dimensions}
Let $(U,d_U)$ be an open,  connected
subset of the Riemann sphere avoiding at least three points and
equipped with a hyperbolic metric.
Let $K\subset U$ be a compact
subset. We denote by $\EE(K,U)$ the space of unramified 
conformal covering maps, $f: \D_f \rr U$,
with the requirement that the covering domain $\D_f\subset K$. 
Denote by $Df:\D_f \rr \RR_+$ the conformal derivative of $f$, see equation
     (\ref{def double lim}),
and by $\|Df\|=\sup_{f^{-1}K} {Df}$
 the maximal value of this derivative over the
set $f^{-1}K$. 
Let $\FF=(f_n) \subset \EE(K,U)$ be a sequence
of such maps.
The intersection 
 \begin{equation}
      J(\FF) = \bigcap_{n\geq 1} f_1^{-1} \circ \cdots \circ f_n^{-1} (U)
 \end{equation}
defines a uniformly hyperbolic
  Julia set for the sequence $\FF$. Let $(\Upsilon,\nu)$ be
a probability space and
let $\omega\in\Upsilon\rr f_\omega\in \EE(K,U)$ be a measurable map.
Suppose that the elements in the
sequence $\FF$ are picked independently, according to the law $\nu$.
Then $J(\FF)$ becomes a random `variable'.
Our main objective is to establish the following

\begin{Theorem}
\label{Main Theorem}
\mbox{}

 I. Suppose that 
${\EEE(\log \,\|Df_\omega\|)}<\infty$.  Then
the Hausdorff dimension of $J(\FF)$ equals almost surely its upper and
lower box dimensions and is 
given by a generalization of Bowen's formula.

 II. Suppose in addition that:
(a) The family of maps $(f_\omega)_{\omega\in \Upsilon}$
  and the probability measure, $\nu$,
depend uniformly real-analytically on complex parameters.
(b)
For any local inverse $f_\omega^{-1}$,
 $\log Df_\omega\circ f_\omega^{-1}$  is uniformly Lipschitz
in parameters and in $\omega\in\Upsilon$.
 (c) The condition number,
         $\|Df_\omega\| \cdot
         \|1/Df_\omega\|$,
is uniformly bounded in parameters and in $\omega\in\Upsilon$.
 Then the almost sure Hausdorff dimension
depends real-analytically
 on the parameters (for more precise definitions
 see section \ref{section real analytic}).
\end{Theorem}

\begin{Example}
\label{example main}
Let $a\in\CC$ and $r\geq 0$ be such that $ |a|+r<\frac{1}{4}$.
Suppose that $c_n\in\CC$, $n\in \NN$ are i.i.d.\ 
 random variables uniformly
distributed in the closed disk $\overline{B}(a,r)$ and that $N_n$, $n\in\NN$ are
 i.i.d.\ random variables distributed according
to a Poisson law of parameter $\lambda \geq 0$.
We consider the sequence of maps, $\FF=(f_n)_{n\in\NNs}$, given by
 \begin{equation}
       f_n(z)= z^{N_n+2} + c_n.
 \end{equation}
% Each map belongs to the space 
% $\EE(\overline{A_{1/2}},A_{1/4+\sigma})$ where
% $A_r=\{r < |z| < 1/r\}$.
An associated `random' Julia set may be defined through
  \begin{equation}
     J(\FF) = \partial \; \{ z \in \CC: f_n \circ \cdots \circ f_1(z)
               \rr \infty \} 
  \end{equation}
As shown in  section  \ref{section real analytic}
the family verifies all conditions
for Theorem  \ref{Main Theorem}, part  I and II.
The random Julia set
therefore has the same almost sure Hausdorff and upper/lower box 
dimension $\dim(J(\FF))=d(a,r,\lambda)$ which in addition
depends real-analytically 
upon $a$,  $r$ and $\lambda$. Note that the
sequence of degrees, $(N_n)_{n\in\NNs}$, almost surely is unbounded
when $\lambda>0$.
\end{Example}

Rufus Bowen, one of the founders of the Thermodynamic Formalism
(henceforth abbreviated  TF), saw
more than twenty years ago
 \cite{Bow79}
a natural connection between 
the geometric properties of a conformal repeller and the
TF for the map(s) generating this repeller.
The Hausdorff dimension \HD($\Lambda$)
of a smooth and compact conformal repeller $(\Lambda,f)$ is precisely
the unique zero $s_\crit$ of a `pressure' function $P(s,\Lambda,f)$ 
having its origin in the TF.
This relationship
 is now known as `Bowen's formula'. The original proof by Bowen
\cite{Bow79} was in the context of Kleinian groups and 
involved a finite Markov partition and
uniformly expanding conformal maps. 
Using TF he constructed a finite Gibbs measure
of zero `conformal pressure' and showed that this measure
is equivalent to the  $s_\crit$-dimensional 
Hausdorff measure of $\Lambda$.
The conclusion then follows.

Bowen's formula apply in many other cases.
For example, when dealing with expanding `Markov maps' 
the Markov partition need not be finite 
and one may eventually have 
a neutral fixed
point in the repeller \cite{Urb96,SSU01}.
One may also relax on smoothness of the maps involved.
 Barreira \cite{Bar96}
and also Gatzouras and Peres \cite{GP97}
were able to demonstrate that Bowen's formula holds
 for   classes of $C^1$ repellers.
A priori, the classical TF does not apply in this setup.
Gatzouras and Peres circumvene the problem by using an
approximation argument and then apply the classical theory.
 Barreira, following the approach of Pesin \cite{Pes88}, 
defines the Hausdorff dimension as a Caratheodory dimension
characteristic. By extending the TF itself Barreira goes closer to
the core of the problem and  may also consider
maps somewhat beyond the $C^1$
case mentioned.  The proofs are, however,  fairly involved
and do not generalize easily neither to a random set-up
nor to a study of parameter-dependency.

In \cite{Rue82}, Ruelle showed that the Hausdorff dimension
of the Julia set of a uniformly  hyperbolic rational map
depends real-analytically
on parameters.  The original approach of Ruelle 
was indirect, using dynamical zeta-functions, \cite{Rue76}.
Other later proofs are based on holomorphic motions,
(see e.g.\ Zinsmeister \cite{Zin99} and references therein).
In either case it  is
difficult to adapt the proofs to a time-dependent and/or random set-up
because the methods do not give sufficiently uniform bounds.
In another context, Furstenberg and Kesten,
\cite{FK60}, had shown, under a condition of log-integrability,
that a random product of matrices has a unique almost
sure characteristic exponent.
 Ruelle, in \cite{Rue79}, required in addition that
the matrices contracted uniformly a positive cone and satisfied
a compactness and continuity condition with respect to the underlying
 propability space.  He showed that under these conditions if the 
family of postive random matrices depends real-analytically
on parameters then so does the almost sure characteristic exponent
of their product.
 He did not, however, allow the probability law to depend on parameters.
 We note here
that if the matrices contract uniformly
a positive cone, the topological conditions in \cite{Rue79}
 may be replaced by the weaker condition of
measurablity + log-integrability.
We also mention 
the  more recent paper, \cite{Rue97}, of Ruelle.
It is in spirit close to \cite{Rue79}
(not so obvious at first sight)
  but provides a more global and far more elegant
 point of view to the question of parameter-dependency.
It has been an invaluable source
 of inspiration to our work.\\

In this article we depart from the traditional path stuck out by TF. In Part
I we  present a proof of Bowen's formula, Theorem
\ref{thm Bowens formula},
 for a  $C^1$ conformal repeller
which bypasses measure theory and most of the TF.
Measure theory can be avoided essentially because $\Lambda$ is compact
and the only element remaining from TF is a family of transfer
operators which encodes geometric informations
into analytic ones.
 Our proof is short and elementary and releases us from
some of the smoothness conditions imposed by TF.
An elementary proof of Bowen's formula should be
 of interest in its own,
at least in the author's opinion.
It generalizes, however, also
to situations where a `standard' approach
either fails or manages only with great difficulties.
We consider  classes of time-dependent 
conformal repellers. By picking
a sequence of maps within a suitable
equi-conformal class one may
study the associated time-dependent repeller. Under the
assumption of uniform equi-expansion and equi-mixing 
and a technical assumption of sub-exponential `growth'
of the involved sequences 
we show, 
Theorem \ref{Thm time dep conf rep},
 that the Hausdorff and Box Dimensions are bounded within the unique zeros
of a lower and an upper conformal pressure. 
Similar results were found by
Barreira \cite[Theorem 2.1 and 3.8]{Bar96}.
When it comes to random conformal repellers, however, the 
approach of Pesin
and Barreira seems difficult to generalize. 
Kifer \cite{Kif96} and later, also
Bogenshutz and Ochs \cite{BO99}, using time-dependent TF
and Martingale arguments,
considered random conformal repellers for certain classes of 
transformations, but under the smoothness restriction imposed
by TF.
In Theorem
\ref{Thm time depend},
 a straight-forward
application of Kingmans sub-ergodic Theorem,
  \cite{King68}, allows us to deal with such cases 
without such restrictions. In addition we obtain very general formulae
for the parameter-dependency  of the Hausdorff dimension.\\

Part II is devoted to Random Julia sets on hyperbolic subsets of
the Riemann sphere.
 Here statements and hypotheses
attain much more elegant forms, cf.\ 
Theorem  \ref{Main Theorem} and
Example \ref{example main} above.  
Straight-forward Koebe estimates enables us to apply Theorem 
\ref{Thm time depend} 
to deduce Theorem \ref{Thm erg analytic}
which in turn yields
Theorem \ref{Main Theorem}, part (I).
The parameter dependency is, however, more subtle.
The central ideas are then the following:
\begin{itemize}
\item[(1)] We introduce a `mirror embedding' of  our hyperbolic subset
    and then a related
   family of transfer operators and 
   cones which a natural (real-)analytic structure. 
\item[(2)] We compute the pressure function as a hyperbolic fixed
   point of a holomorphic map acting upon  the cone-family. 
   When the family of maps depends real-analytically   on parameters,
   then the real-analytically dependency of the dimensions,
    Theorem \ref{part I},
    follows from an implicit function theorem.
\item[(3)] The above mentioned fixed point is hyperbolic.
   This implies an exponential decay of the fixed point with `time'
    and allow us to treat a real-analytic parameter dependency
   with respect to  the underlying probability law. 
   This concludes the proof of Theorem \ref{Main Theorem}.
\end{itemize}

\newpage

\section{Part I: $C^1$ conformal repellers and Bowen's formula}
Let $(\Lambda,d)$ be a non-empty compact metric space
without isolated points
and let $f:\Lambda \rr \Lambda$ be a continuous surjective map.
Throughout Part I we will write interchangeably
 $f_x$ or $f(x)$ for a
map $f$ applied to a point $x$.
We say that $f$ is {\em $C^1$ conformal} at $x\in \Lambda$ iff
the following double limit exists~:
  \begin{equation}
     Df_x = \lim_{u \neq v\rr x} \frac{d(f_u,f_v)}{d(u,v)}.
     \label{def double lim}
  \end{equation}
The limit is called the conformal derivative of $f$ at $x$.
The map $f$ is said to be $C^1$ conformal on
$\Lambda$ if it is so at every
point of $\Lambda$.
A point $x\in\Lambda$ is said to be {\em critical} if 
$ a_x=0$.

The product $Df^{(n)}_x=Df_{f^{n-1}x} \cdots Df_x$ along the orbit of $x$ is
the conformal
derivative for the $n$'th iterate of $f$.
The map is said to be uniformly expanding if
  there are constants $C>0$, $\beta\!>\!1$ for which
        $Df^{(n)}_x \geq C \beta^n$
         for all $x\in\Lambda$ and $n\in\NN$.
We say that $(\Lambda,f)$ is a  $C^1$ conformal repeller
if
\begin{enumerate} \isep
\item[(C1)] $f$ is {$C^1$ conformal} 
       on $\Lambda$.
\item[(C2)]  $f$ is {uniformly expanding}.
\item[(C3)] $f$ is an open mapping.\\
\end{enumerate}

For $s\in\RR$ we define
 the dynamical pressure
 of the $s$-th power of the conformal derivative
by the formula:
\begin{equation}
   P(s,\Lambda,f) = \liminf_n\  \frac{1}{n} \log \
           \sup_{y\in\Lambda} \
	   \sum_{x\in\Lambda:f^n_x=y}
	       \left( {Df^{(n)}_{x}} \right)^{-s}.
	       \label{conf deriv}
\end{equation}

We then have the following
\begin{Theorem}\mbox{\rm\bf (Bowen's formula)}
    \cite{Bow79,Rue82,Fal89,Bar96,GP97} Let $(\Lambda,f)$ be 
a  $C^1$ conformal repeller.
Then, the Hausdorff dimension of
 $\Lambda$ coincides with its upper and lower box dimensions and 
is given as the unique zero of the pressure function $P(s,\Lambda,f)$.\\
\label{thm Bowens formula}
\end{Theorem}

For clarity of the proof we will here impose the additional
assumption of strong mixing.
 We have delegated to
Appendix \ref{app removing} a sketch of how to remove this restriction.
We have chosen to do so because (1) the proof is really much more
elegant  and (2) there seems to be no
natural generalisation when dealing with the time-dependent case
(apart from trivialities).

More precisely, to any given
$\delta>0$
we assume that
there is an integer
$n_0=n_0(\delta)<\infty$
(denoted the $\delta$-covering
time for the repeller)
 such that for every $x\in\Lambda$ :
\begin{enumerate}
\item[(C4)] \ \ 
         \begin{equation}
	 f^{n_0}B(x,\delta)=\Lambda .
         \label{mixing}
	 \end{equation}
\end{enumerate}
For the rest of this section, $(\Lambda,f)$ will be assumed to be
a strongly mixing $C^1$ conformal repeller, thus verifying 
(C1)-(C4).\\

Recall that
a countable family 
$\{U_n\}_{n\in \NN}$
of open sets 
is a
$\delta$-cover($\Lambda$) if $\diam\;  U_n <\delta$ for all $n$ and
their union contains (here equals) $\Lambda$. For $s\geq 0$ we set
\[
   M_\delta(s,\Lambda) =
    \inf\ \left\{ \sum_n (\diam\; U_n)^s :
    \{U_n\}_{n\in\NN} \mbox{\ is a \ } \delta\!-\!\mbox{cover}(\Lambda) 
     \right\} \in [0,+\infty]
\]
Then $M(s,\Lambda)=\lim_{\delta\rr 0} 
 M_\delta(s,\Lambda) \in [0,+\infty]$ exists and is called
the $s$-dimensional Hausdorff measure of $\Lambda$.
The Hausdorff dimension
 is the unique critical  value $s_\crit=\HD \Lambda\in[0,\infty]$
such that $M(s,\Lambda)=0$ for $s>s_\crit$ and
$M(s,\Lambda)=\infty$ for $s<s_\crit$. The Hausdorff measure is said
to be finite if $0<M(s_\crit,\Lambda)<\infty$.

Alternatively we may replace the condition on the covering sets
by considering finite covers by open balls $B(x,\delta)$  of fixed radii
 $\delta>0$.
Then the limit as $\delta\rr 0$ of $M_\delta(s,\Lambda)$
need not exist so we replace it by taking lim sup and lim inf.
We then obtain the upper, respectively the lower $s$-dimensional
Box `measure'. The upper and lower Box Dimensions, $\boxsup \Lambda$
and $\boxinf \Lambda$, are the corresponding critical values.
It is immediate that

\[
 0\leq   \dimH \Lambda \leq \boxinf \Lambda \leq \boxsup \Lambda \leq +\infty
\]

\begin{Remarks}\mbox{}\\[-5mm]
\begin{enumerate}
  \item Let $J(f)$ denote the
      Julia set  of a uniformly hyperbolic rational map
       $f$ of the Riemann sphere. There is an open
       (hyperbolic) neighborhood $U$ of $J(f)$  such that
       $V=f^{-1}U$ is compactly contained in $U$ and 
    such that $f$ has no critical points
       in $V$.  Writing $d$ for the hyperbolic metric on $U$
        one verifies that $(J(f),f)$
       is a $C^1$ conformal repeller. 
\item
  Let $X$ be a $C^1$ Riemannian manifold without boundaries and let
  $f:X\rr X$ be a $C^1$ map. It is an exercise in Riemannian 
  geometry to see that
  $f$ is uniformly conformal at $x\in X$
  iff $f_{*x}:T_xX \rr T_{fx}X$ is a conformal map of tangent spaces and in
  that case, $Df_x = \|f_{*x}\|$. 
  When dim $X < \infty$
  condition (C3) follows from (C1)-(C2).
  We note also that
  being $C^1$ 
(the double limit in 
     equation \ref{def double lim})
rather than just differentiable is important.
\end{enumerate}
\end{Remarks}

\subsection{Geometric bounds}
We will first establish sub-exponential
geometric bounds for iterates
of the map $f$.
In the following we say that
a sequence $(b_n)_{n\in\NN}$ of positive real numbers
is sub-exponential or 
 of sub-exponential growth if $\lim_n \sqrt[n]{b_n}=1$.
For notational convenience we will also assume that
$Df_x\geq \beta>1$ for all $x\in\Lambda$. This may always be
achieved in the present set-up by considering a
high enough iterate of the map $f$, possibly redefining $\beta$.

Define the {\em divided difference},
\begin{equation}
  f[u,v]=\left\{
         \begin{array}{ll}
             \frac{d(f_u,f_v)}{d(u,v)} \ \ \ & u\neq v \in \Lambda,\\
              Df_u   & u=v\in \Lambda.
         \end{array} \right.
   \label{divided diff}
\end{equation}
Our hypothesis on $f$ implies that $f[\cdot,\cdot]$ is continuous
 on the compact set $\Lambda\times \Lambda$ 
and not smaller than $\beta>1$ on the diagonal of the product set.
 We let
 $\|Df\|=\sup_{u\in \Lambda} Df_u < +\infty$ denote the maximal
conformal derivative on the repeller.

Choose $1<\lambda_0<\beta$. Uniform continuity and openness of the map $f$
show that we may find $\delta_f>0$ and $\lambda_1<+\infty$
such that

 \begin{enumerate}\isep
  \item[(C2')] \ \ \
              $\lambda_0 \leq f[u,v]\leq \lambda_1 $ \ \ \
         whenever\ \  \
        $d(u,v)<\delta_f$.
  \item[(C3')]\ \ \
         $ B(f_x,\delta_f) \subset f B(x,\delta_f)$, for all $x\in\Lambda$.
 \end{enumerate}

The constant $\delta_f$ gives a scale below which the
map $f$ is injective, uniformly expanding and (locally) onto.
In the following we will assume that values of $\delta_f>0$,
$\lambda_0>1$ and $\lambda_1<+\infty$ have been found so as to
satisfy conditions (C2') and (C3').\\

We define the {\em distortion} of $f$ at $x\in\Lambda$
and for $r>0$  as follows:
\begin{equation}
      \epsilon_f(x,r)=
           \sup \{
             \ \log \; 
	     \frac{f[u_1,u_2]}
	          {f[u_3,u_4]} : {\rm all} \ 
           u_i\in B(x,\delta_f)\cap f^{-1}B(f_x,r)
                  \}.
      \label{def distortion}
   \end{equation}
 This quantity tends to zero as $r \rr 0^+$ uniformly in $x\in\Lambda$
(same compactness and continuity as before).
Thus,
   \[
        \epsilon(r) = \sup_{x \in \Lambda} \epsilon_f(x,r)
   \]
tends to zero as $r \rr 0^+$.
When $x\in\Lambda$ and the $u_i$'s are as in (\ref{def distortion})
then also:
\begin{equation}
   \left| \log\ \frac{f[u_1,u_2]}{Df_{u_3}} \right| \leq 
              \epsilon(r) 
    \mbox{\ \ \ \ and\ \ \ \ }
    \label{first dist bound}
   \left| \log\ \frac{Df_{u_1}}{ Df_{u_2}} \right| \leq 
              \epsilon(r) .
     \label{second dist bound}
\end{equation}

\mbox{}\\

For $n\geq \NN$ we 
define 
the n-th `Bowen ball' around $x\in\Lambda$,
\[
    B_n(x) \equiv B_{n}(x,\delta_f,f)=
      \{ u\in\Lambda: d(f^k_x,f^k_u) < \delta_f, \ 0\leq k \leq n\}
            .
\]
We say that $u$ is $n$-{\em close} to 
$x\in\Lambda$ if $u\in B_n(x)$.
The Bowen balls act as `reference' balls, getting
uniformly smaller with increasing $n$.
In particular,
$\diam\;B_n(x) \leq 2\; \delta_f \; \lambda_0^{-n} $, i.e.\
tends to zero exponentially fast with $n$.
We also see that for each $x\in\Lambda$ and $n\geq 0$ the map,
  \[
    f\;: \; B_{n+1}(x) \rr B_n(f_x),
 \]
is a uniformly expanding homeomorphism. 

Expansiveness of the map $f$ means that closeby points
may follow very different future trajectories.  Our assumptions
assure, however,  that closeby points have 
very similar backwards histories. 
The following two Lemmas emphasize this point~:

\begin{Lemma}{\bf [Pairing]}
For each $y,w\in\Lambda$ with $d(y,w)\leq \delta_f$ and for
every $n\in \NN$  the sets $f^{-n}\{y\}$ and $f^{-n}\{z\}$
may be paired uniquely into pairs of
$n$-close points.
\end{Lemma}

Proof: Take $x\in f^{-n}\{y\}$. 
The map 
$f^n: B_n(x) \rr B_0(f^n_x)=B(y,\delta_f)$
is a homeomorphism. Thus there is a unique point 
$u\in f^{-n}\{z\} \cap B_n(x)$. By construction,
$x\in B_n(u)$ iff $u\in B_n(x)$. Therefore
$x\in f^{-n}\{y\} \cap B_n(u)$ is the unique pre-image of
 $y$ in the n-th Bowen ball
around $u$ and we obtain the desired pairing.
 \Halmos\\

\begin{Lemma}{\bf [Sub-exponential Distortion]}
\label{lemma distortion}
There is a sub-exponential sequence,
 $(c_n)_{n\in\NN}$,
 such that for any two points $z,u$
 which are $n$-close to  $x\in\Lambda$ ($x \neq u$)
  \[
     \frac{1}{c_n} \leq 
               \frac
               {d(f^n_u,f^n_x)} 
	       {d(u,x) \; Df^{(n)}_z} 
	       \leq c_n
      \mbox{\ \ \ \ and \ \ \ \ }
     \frac{1}{c_n} \leq 
               \frac
               {Df^{(n)}_x}
	       {Df^{(n)}_z}
	       \leq c_n
  \]
\end{Lemma}

Proof:  For all $1\leq k \leq n$ we have that
$f^{k}_u\in B_{n-k}(f^k_x)$.
Therefore, $d(f^k_u,f^k_x)< \delta_f \lambda_0^{k-n}$ and
the distortion bound
(\ref{first dist bound}) implies that
\[
   |\log 
               \frac
               {d(f^n_u,f^n_x)} 
	       {d(u,x) \; Df^{(n)}_z} | \leq 
               \epsilon(\delta_f) + \epsilon(\delta_f \lambda_0^{-1}) +
              \cdots+ 
	       \epsilon(\delta_f \lambda_0^{1-n})
               \equiv \log  c_n.
   \]
Since $\lim_{r\rr 0}\epsilon(r)= 0$ 
it follows that $\frac1n \log c_n \rr 0$, whence 
that the sequence $(c_n)_{n\in\NN}$ is of sub-exponential growth.
This yields the first inequality and
 the second is proved
e.g.\ by taking the limit $u\rr x$.
\Halmos\\

\begin{Remarks}
 When $K=\int_0^{\lambda_0 \delta_f} \epsilon(t)/t \; dt < +\infty$
 one verifies that the distortion stays uniformly bounded,
i.e.\ that $  c_n \leq K/(\lambda_0-1) < \infty$ uniformly in $n$. This is
 the case, e.g.\ when $\epsilon$ is H\"older continuous at zero.
\label{remark uniform cn}
\end{Remarks}

\subsection{Transfer operators}
Let $\calM(\Lambda)$ denote the Banach space of bounded real valued 
functions on $\Lambda$ equipped with the sup-norm. We denote by
 $\chiU$ the characteristic function of a subset $U\subset \Lambda$
and we write 
$\bfone=\chi_\Lambda$ 
for the constant function $\bfone(x)=1$, $\forall x\in\Lambda$.
For $\phi\in \calM(\Lambda)$
 and $s\geq 0$ we define the positive linear transfer\footnote{
      The `transfer'-terminology, inherited from statistical mechanics,
      refers here to the `transfer' of the encoded geometric information at
      a small
      scale to a larger scale, using the dynamics of the map, $f$.}
operator
\[
   (L_s \phi)_y \equiv (L_{s,f} \phi)_y \equiv \sum_{x\in \Lambda: f_x=y} 
        \left( {Df_x} \right)^{-s}
        \phi_x,\ \ \  y\in\Lambda.
        \label{def trop1}
\]
Since $\Lambda$ has a finite $\delta_f$-cover and $Df$
is bounded these operators are necessarily bounded.
The $n$'th iterate of the operator $L_s$ is given by
   \[
   (L_s^n \phi)_y = \sum_{x\in \Lambda: f^n_x=y} 
        \left( {Df^{n}_x} \right)^{-s}
        \phi_x.
\]
It is of importance to obtain bounds for the action of $L_s$ upon the
constant function. More precisely,
for  $s\geq 0$ we denote
  \begin{equation}
      M_n(s)\equiv \sup_{y\in\Lambda} L^n_s \bfone(y)
      \ \ \ \ \mbox{and} \ \ 
      \ \ m_n(s)\equiv \inf_{y\in\Lambda} L_s^n \bfone (y) .
      \label{lower upper}
   \end{equation}
We then define the lower, respectively the upper pressure through
\[
  -\infty \leq {\underline{P}(s)} \equiv \liminf_n \frac1n \log m_n(s) \ \ \ 
     \leq \ \ \ 
  {\overline{P}(s)} \equiv \limsup_n \frac1n \log M_n(s) \leq +\infty.
\]
\begin{Lemma} {\bf [Operator bounds]}
  For each $s\geq 0$, the upper and lower pressures agree and are finite.
  We write $P(s)\equiv \UP(s)=\OP(s)\in \RR$.
 for the common value.
  The function  $P(s)$ is continuous, strictly decreasing
  and  has a unique zero, $s_\crit\geq 0$.
  \label{Operator bounds}
\end{Lemma}
Proof: Fix $s\geq 0$.
Since the operator is positive, the sequences $M_n=M_n(s)$
and $m_n=m_n(s)$, ${n\in\NN}$ are 
sub-multiplicative and super-multiplicative, respectively.
Thus,
\begin{equation}
 m_k m_{n-k} \leq m_n \leq M_n \leq M_{k}  M_{n-k}
                , \ \ \ \forall 0\leq k\leq n .
   \label{sub super mult}	
\end{equation}
This implies convergence of both $\sqrt[n]{M_n}$ and
$\sqrt[n]{m_n}$,
the limit of the former sequence being the spectral
 radius of $L_s$ acting upon $\calM(\Lambda)$.
Let us sketch a standard proof for the first sequence:
Fixing $k\geq 1$ we write
$n=pk+r$ with $0\leq r<k$.
Since $k$ is fixed, $\lim \sup_n \max_{0\leq r <k}
\sqrt[n]{M_r} = 1$. But then
$\lim \sup_n \sqrt[n]{M_n}=
  \lim \sup_p 
    \sqrt[pk]{M_{pk}} \leq
\sqrt[k]{M_k}$. Taking lim inf (with respect to $k$) on the right hand side
we conclude that the limit exists. A similar proof works for the 
sequence $(m_n)_{n\in\NNs}$.
Both limits are non-zero $(\geq m_1>0)$ and finite
 $(\leq M_1<\infty)$. 
We need to show that the ratio
$M_n/m_n$ is of sub-exponential growth.

 Consider $w,z\in\Lambda$ with $d(w,z)< \delta_f$ and $n>0$.
The Pairing Lemma shows that we may pair
the pre-images $f^{-n}\{w\}$ and $f^{-n}\{z\}$ 
into pairs of $n$-close points, 
say $(w_\alpha,z_\alpha)_{\alpha\in I_n}$ over
a finite index set  $I_n$, possibly depending on the pair $(w,z)$.
Applying the second distortion bound in
Lemma \ref{lemma distortion}
to each pair yields
\begin{equation}
    L^n_s \; \bfone (z) \geq \left(\frac{1}{c_n}\right)^s L^n_s \; \bfone (w).
    \label{w z}
\end{equation}

Choose $w\in\Lambda$ such that
\ $L^n_s \; \bfone(w)\geq M_n/2$.
Given an arbitrary $y\in\Lambda$  our strong mixing assumption 
(C4)
implies that
the set $B(w,\delta_f)\cap f^{-n_0}\{y\}$
contains at least one point.
 Using (\ref{w z})
 we obtain
\[
L^{n+n_0}_s \; \bfone (y) =
           \sum_{z:f^{n_0}_z=y} \left( {Df^{n_0}_z} \right)^{-s} 
            L^n_s \; \bfone(z) 
               \geq 
      \left( {\|Df\|^{n_0} c_n} \right)^{-s} 
          \frac{M_n}{2} .
\]
Thus, 
\begin{equation}
    m_{n+n_0} \geq 
 (\|Df\|^{n_0} c_n)^{-s} M_n/2
  \label{mM bounds}
\end{equation}
 and since $c_n$ is of sub-exponential
growth then so is $M_n/m_{n+n_0}$ and 
therefore also $M_{n+n_0}/m_{n+n_0}\leq M_{n_0}M_n/m_{n+n_0} $.\\

The functions $s\log \beta+P(s)$ 
and $s\log \|Df\|+P(s)$ are non-increasing and non-decreasing,
respectively. Also $0\!\leq\! P(0)\!<\! +\infty$.
  It follows that $s\mapsto P(s)$ is continuous
 and that $P$ has a unique zero $s_\crit\geq 0$. \Halmos

\begin{Remarks}
  Super- and sub-multiplicativity (\ref{sub super mult})
    imply the 
  bounds\footnote{Such bounds are useful in applications as they yield
        computable rigorous bounds for the dimensions.}
  \[
           m_n(s)  \leq e^{nP(s)} 
      \leq  M_n(s), \ \ n\in\NN.
  \]
Clearly, if the distortion $c_n$ is uniformly bounded then so is the ratio, 
	$M_n/m_n \leq K(s) < \infty$.
  \label{remark uniform bn}
\end{Remarks}

To prove Theorem \ref{thm Bowens formula}
 it suffices to show that 
$s_\crit\leq \HD(\Lambda)$ and $\boxsup(\Lambda)\leq s_\crit$.
%In order to do so our main
%ingredient will be to consider the action
%of a transfer operator on a characteristic function
%on a subset $U$. When the diameter of $U$ is small and $k$ is not too large 
%so that distortion is under control we establish a relationship
%of the form ($\sim$ means here within sub-exponential bounds):
%\begin{equation}
%      L_s^k \chiU \sim \left( \frac{\diam\; U}{\diam \; f^k U} \right)^s
%                 \chi_{f^k U}.
%\end{equation}
%The characteristic function on the right hand side is on a `large' subset

\subsection{\HD($\Lambda$)$\geq s_\crit$}
\label{lower bounds}
Let $U\subset \Lambda$ be an open non-empty subset of diameter
not exceeding $\delta_f$.
We will iterate $U$ by $f$ until the size of $f^kU$ becomes `large'
compared to $\delta_f$. 
As long as $f^k$ stays injective on $U$ the set
$\{ z\in U: f^k_z=y\}$ contains at most one element for any $y\in\Lambda$.
Therefore, for such $k$-values 
  \begin{equation}
     (L^k_s \chiU)(y) \leq \sup_{z\in U} 
            \left( {Df^{k}_z} \right)^{-s}
   \ , \ \ \forall \ y\in\Lambda.
      \label{single bound}
  \end{equation}
Choose $x=x(U)\in U$ and let
$k=k(U)\geq 0$ be the largest positive integer
for which $U\subset B_k(x)$. In other words~:
\begin{enumerate}\isep 
\item[(a)] $d(f^l_x,f^l_u) <  \delta_f$ for $0\leq l \leq k$ 
             and all  $u\in U$,
\item[(b)] $d(f^{k+1}_x,f^{k+1}_u) \geq \delta_f$ for some $u\in U$.
\end{enumerate}
$k(U)$ is finite because $U\cap \Lambda$ contains at 
least two distinct points
which are going to be separated when iterating.
Because of (a) 
 $f^k$ is injective on $U$ so that
      (\ref{single bound}) applies.
On the other hand, (a) and (b) implies that there is $u\in U$ for which 
$\delta_f \leq d(f^{k+1}_x,f^{k+1}_u) 
\leq \lambda_1(f) d(f^k_u,f^k_x) $ where 
$\lambda_1(f)$ was the maximal dilation of $f$ on $\delta_f$-separated points.
Our sub-exponential distortion estimate shows that for any $z\in U$,
  \[
   \frac{\delta_f/\lambda_1(f)}{ 
   \diam \;U }\; \frac{1}{Df^{k}_z}  \leq 
   \frac{d(f^k_u,f^k_x)}{d(u,x)}\;\frac{1}{Df^{k}_z} \leq c_k.
  \]
Inserting this in (\ref{single bound}) and using the 
definition of $m_n(s)$   we see that for any $y\in\Lambda$,
  \begin{equation}
   (L^k_s \chiU)(y) \leq \;
          (\diam\; U)^s \; (\frac{\lambda_1(f) c_k}{\delta_f})^s\; \bfone \leq 
             (\diam\; U)^s 
            \left[ (\frac{\lambda_1(f) c_k}{\delta_f})^s\;
                  \frac{1}{m_k(s)}\right] \; L^k_s \bfone .
	     \label{Lk bound}
  \end{equation}
Choosing now $0<s<s_\crit$, the sequence $m_k(s)$ tends exponentially
fast to infinity
[when $s_\crit=0$ there is nothing to show].
Since the sequence  $\left((c_k)^s\right)_{k\in\NNs}$ 
is  sub-exponential the factor in square-brackets
is uniformly bounded in $k$, say by $\gamma_1(s)<\infty$ (independent of $U$).
Positivity of the operator implies that for any $n \geq k(U)$
  \[
     L^n_s \chiU \leq \gamma_1(s) \ 
          (\diam\; U)^s \ 
           L^n_s \bfone .
   \]

If $(U_\alpha)_{\alpha\in\NN}$
 is an open $\delta_f$-cover of the compact set $\Lambda$ then
it has a finite sub-cover, say $\Lambda \subset U_{\alpha_1} \cup \ldots
       \cup U_{\alpha_m}$. 
Taking now $n = \max\{k(U_{\alpha_1}),\ldots,k(U_{\alpha_m})\}$ 
 we obtain
  \begin{equation}
     L^n_s \bfone \leq \sum_{i=1}^m L^n_s \chi_{U_{\alpha_i}} 
             \leq \gamma_1(s)\  
              \sum_{i=1}^m (\diam\; U_{\alpha_i})^s\
                L^n_s \bfone 
                  \leq \gamma_1(s) \ 
        \sum_\alpha (\diam\; U_{\alpha})^s\
		  L^n_s \bfone .
		  \label{lower bound}
  \end{equation}

This equation shows that $\sum_\alpha (\diam\; U_\alpha)^s$ is 
bounded uniformly from below by $1/\gamma_1(s)>0$.
 The Hausdorff dimension of $\Lambda$
is then not smaller than $s$, whence not smaller than $s_\crit$.

\subsection{$\boxsup \Lambda \leq s_\crit$}
\label{dimB upper}
\begin{figure}
\begin{center}
\epsfig{figure=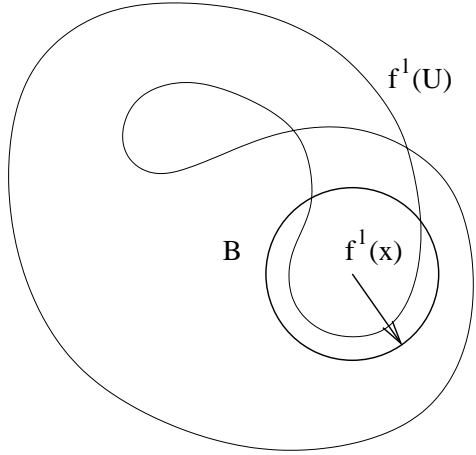,height=6cm}
\label{box-dim}
\end{center}
\caption{An iterate $f^l(U)$ which covers $B=B(f^l(x),\delta_f)$ but not
in the `right' way.}
\end{figure}

Fix $0<r<\delta_f$
and 
let $x\in\Lambda$.
 This time we wish
to iterate a ball  $U={B(x,r)}$ until it has a `large' interior and
contains a ball of size $\delta_f$. 
This may, however, not be good enough 
(cf.\ Figure \ref{box-dim}).
We also need to control the distortion.
Again  these two
goals combine nicely when considering 
the sequence of Bowen balls,
$B_k\equiv B_k(x)$, $k\geq 0$. 
It forms a sequence of neighborhoods of
$x$, shrinking to $\{x\}$.
Hence, there is a smallest integer $k=k(x,r)\geq 1$
such that $B_{k} \subset U$.
$f^k$ maps $B_k$ homeomorphically onto
 $B_0(f^k_x)=B(f^k_x,\delta_f)$ and
positivity of $L_s$ shows that
\[
   L_s^k\chi_{U}\geq L_s^k \chi_{B_{k}} \geq
     \inf_{z\in B_{k}} 
     \left( {Df_z^{k}} \right)^{-s} \chi_{B(f^k_x,\delta_f)} .
\]

By assumption $B_{k-1}\not\subset U$ 
so there must be a point
$y\in B_{k-1}$ with $d(y,x)\geq r$.
%But then 
%$d(f^k_u,f^k_x)\leq \lambda_1 d(f^{k-1}_u,f^{k-1}_x)
 %\leq \lambda_1 \delta$.
As $y$ is 
$(k-1)$-close to $x$ 
our distortion estimate shows that for any 
$z\in B_k(x)\subset B_{k-1}(x)$,
  \[
   \frac{\delta_f}{ r} \frac{\|Df\|}{Df^{k}(z)} >
   \frac{d(f^{k-1}_y,f^{k-1}_x)}{d(y,x)}
       \; \frac{1}{Df^{{k-1}}(z)} \geq \frac{1}{c_{k-1}}.
    \label{upp conf bd}
  \]

Therefore,
\[
L_s^k \chi_{U} \geq 
         r^s (\delta_f c_{k-1} \|Df\|)^{-s} \chi_{B(f^k_x,\delta_f)} .
\]
If we iterate another $n_0=n_0(\delta_f)$ times then $f^{n_0} B(f^k_x,\delta_f)$
 covers all of $\Lambda$  due to mixing and using the definition
of $M_n(s)$,

\[
  L_s^{k+n_0} \chi_{U} \geq 
      r^s (\delta_f c_{k-1} \|Df\|^{1+n_0})^{-s} \bfone \geq
   (4r)^s \left[ \frac
         {(4\|Df\|^{1+n_0}\delta_f c_{k-1} )^{-s}}
          {M_{k+n_0}(s)} \right]
    L_s^{k+n_0} \bfone.
\]

When $s>s_\crit$, $M_{k+n_0}(s)$ tends expontially fast to zero.²
As the rest is sub-exponential,
the quantity in the square brackets is
uniformly bounded from below by some $\gamma_2(s)>0$.
Using the positivity
of the operator we see that
\begin{equation}
       L_s^{n} \chi_{U} \geq \gamma_2(s) (4r)^s  L_s^n \bfone,
     \label{box lower}
\end{equation}
whenever $n\geq k(x,r)+n_0$.

Now, let  $x_1,\ldots,x_{N}$ be a finite maximal $2r$ separated set
 in $\Lambda$.
Thus, the balls $\{B(x_i,2r)\}_{i=1,\ldots,N}$ cover $\Lambda$
whereas the balls $\{B(x_i,r)\}_{i=1,\ldots,N}$ are mutually disjoint.
For $n \geq \max_i k(x_i,r)+n_0$,
   \[
      L_s^n \bfone \geq
       \sum_i L^{n}_s \chi_{B(x_i,r)} \geq \gamma_2(s)\;  N (4r)^s
               L_s^n \bfone.
  \]
 We have deduced the bound
  \[
              \sum_{i=1}^N (\diam\; B(x_i,2r))^s \leq
              {1/\gamma_2(s)}
  \]
which shows that
$\boxsup \Lambda$ does not exceed $s$, whence not $s_\crit$.
We have proven Theorem
\ref{thm Bowens formula}
 in the case of a strongly mixing repeller and refer 
to Appendix  
 \ref{app removing} for the extension to the general case.
\Halmos

\begin{Corollary}
If $\int_0^{\lambda_0 \delta_f} \epsilon(t)/t \; dt <+\infty$
and the repeller is strongly mixing
(cf.\ Remark \ref{non finite measure})
 then the Hausdorff measure is finite 
and comprised between $1/\gamma_1(s_\crit)>0$ and 
$1/\gamma_2(s_\crit)< + \infty$.
\end{Corollary}
Proof: The hypothesis implies that for fixed $s$
the sequences $(c_n(s))_n$ and $M_n(s)/m_n(s)$
 in the sub-exponential distortion
and operator bounds, respectively, are both uniformly bounded in $n$
(Remarks \ref{remark uniform cn} and
\ref{remark uniform bn}). 
All the (finite) estimates 
may then be carried out at $s=s_\crit$ and the conclusion follows.
(Note that no measure theory was used to reach this conclusion).
 \Halmos\\

%\begin{Remarks}
% There is no need here to put more structure upon the space $\calM(\Lambda)$.
% This would be reasonable, however, when studying parameter dependency
% (more than mere continuity) of the Hausdorff dimension.
% In that case one may consider the operator acting upon
% a  suitable subspace of  $C(\Lambda)$ and hope for 
% the presence of a `spectral gap'. The largest (isolated) eigenvalue
% then depends on parameters in the same way as the operator.
% It is convenient also in this case to avoid TF in order to
% preserve maximal smoothness of the operator. As alternatives
% one may use cone-contractions as in \cite{FS79,Liv95}.
%\end{Remarks}

\section{Time dependent conformal repellers}
Let $(K,d)$ denote a 
 complete metric space without isolated points
and let $\Delta>0$ be  such that 
$K$ is covered by a finite number, say
  $N_\Delta$ balls
of size $\Delta$.
To avoid certain pathologies we will also assume 
that 
 $(K,d)$ is  {\em $\Delta$-homogeneous}, i.e.\ that 
  there is a constant  $0<\delta<\Delta$
 such that for any $y\in K$,
  \begin{equation}
 B(y,\Delta)\setminus B(y,\delta)\neq \emptyset. 
   \label{homogeneous}
  \end{equation}
For example, if $K$ is connected or consists of a finite number
of connected components then $K$ is
% {\em per se}
$\Delta$-homogeneous.\\

Let $\beta>1$ and let $\epsilon: [0,\Delta] \rr [0,+\infty[$ be
an $\epsilon$-function, i.e.\
a continuous function with $\epsilon(0)=0$.
Consider the class 
$\EE=\EE(\Delta,\beta,\epsilon)$
of maps $f$ where
\[
      f: \Omega_f \rr K
\]
is a $C^1$-conformal
unramified covering map of finite maximal degree,
$d^o_{\max}(f)= \max_{y\in K} \mbox{deg}(f;y)\in\NN$,
from a non-empty (not necessarily connected) 
domain, $\Omega_f\subset K$, onto $K$, subject to
the following `equi-uniform' requirements: 
There are constants $\delta(f)>0$, $\lambda_1(f)<+\infty$  and a 
 function  $\delta_f :x\in \Omega_f \rr [\delta(f),\Delta] \subset \RR_+$
such that :
\begin{Assumption}
\label{T0-T3}
\mbox{}

\begin{enumerate}
  \item[(T0)]\ \ \ For all distinct $x,x'\in f^{-1}\{y\}$
     (with  $y\in K$) \ \
      the balls $B(x,\delta_f(x))$ and $B(x',\delta_f(x'))$
      are disjoint,\ \ \ (local injectivity),
  \item[(T1)]\ \ \ For all $x \in \Omega_f$  :  \ 
       $B(f_x,\Delta)\ \subset f ( B(x,\delta_f(x)) \cap \Omega_f))$,\ \ \
      (openness),
   \item[(T2)]\ \ \
      $\beta \leq f[u,x] \leq \lambda_1(f)$, \ \ $\forall\; u,x\in \Omega_f$ :  
           $d(u,x)<\delta_f(x)$, 
  (dilation) and
\item[(T3)] \ \ \ For all $x\in\Omega_f$ : \
        $\epsilon_f(x,r) \leq \epsilon(r)$,
       $\forall\  0< r\leq \Delta$,\ \ \
(distortion).\\
\end{enumerate}
\end{Assumption}

Here, $f[\cdot,\cdot]$ is the divided
difference from equation (\ref{divided diff})
and the distortion, a restricted version of
equation (\ref{def distortion}), for $x\in\Omega_f$, $r>0$, is given by
\begin{equation}
      \epsilon_f(x,r)=
           \sup \{
             \ |\log \; 
	     \frac{f[u_1,x]}
	          {Df(u_2)}| :  \
           u_1,u_2\in B(x,\delta_f(x))\cap f^{-1}B(f_x,r)
                  \}.
      \label{def new distortion}
   \end{equation}
We tacitly understand by writing $f^{-1}(y)$ that we are looking at
the pre-images of $y\in K$ within $\Omega_f$, i.e.\ where the map is defined.
We also write $\|Df\|$ for the supremum of the conformal derivative
of $f$ over its domain of definition, $\Omega_f$.  By (T2) and setting $u=x$,
 we also see that $\beta\leq \|Df\|\leq \lambda_1(f)$.

When  $f\in \EE(\Delta,\beta,\epsilon)$ and $f(x)=y\in K$ then by
$\Delta$-homogeneity, (\ref{homogeneous}),
 and property (T1), there must be
$u\in B(x,\delta_f(x))$ with $f(u)\in B(y,\Delta)\setminus B(x,\delta)$
($\delta$ as in the above definition).
By the definition, equation
(\ref{def new distortion}), of the distortion it follows that
 \begin{equation}
 0< \kappa \equiv \delta e^{-\epsilon(\Delta)}  \leq
     \delta_f(x) Df(x) , \ \ \ \forall x\in  \Omega_f .
   \label{def k}
    \end{equation}

In the following let 
$\FF=(f_k)_{k\in\NN} \subset \EE(\Delta,\beta,\epsilon)$
 be a fixed sequence of such mappings and let us fix
 $\delta_{f_k}(x)$,
$\delta_{f_k}=\inf_{x\in\Omega_{f_k}}\delta_{f_k}(x)>0$
 and $\lambda_1(f_k)$ 
 so as to satisfy conditions (T0)-(T3).
For all $n\geq 0$ define~:
\[
    \Omega_n(\FF)= f_1^{-1} \circ \cdots \circ f_n^{-1}(K),
\]
and then 
\[
    \Lambda(\FF)=\bigcap_{n\geq 1}\; \Omega_n(\FF).
\]
% For all $t\geq 0$ define~:
% \begin{equation}
   % \Lambda_t = \Lambda_t(\FF) \equiv \bigcap_{n\geq 1} 
        % f_{t+1}^{-1} \circ \cdots \circ f_{t+n}^{-1}(K) .
% \end{equation}
Letting $\sigma(\FF)=(f_{k+1})_{k\in\NNs}$ denote the shift
of the sequence we set
$\Lambda_t=\Lambda(\sigma^t( \FF))$, $t\geq 0$.
$K$ was assumed complete (though not necessarily compact)
and each $\delta(f_k)$ is strictly positive. It follows then that
each $\Lambda_t$ is closed, whence complete. Each $\Lambda_t$ also
 has finite open covers
of arbitrarily small diameters (obtained by pulling back
 a finite $\Delta$-cover of $K$),
 whence each $\Lambda_t$ is compact and non-empty.
Also $f_t(\Lambda_{t-1})=\Lambda_{t}$ and we have obtained
a time-dependent sequence of compact conformal repellers~:
\[
   \Lambda_{0}
   \stackrel{f_{1}}{\longrightarrow}
   \Lambda_{1}
   \stackrel{f_{2}}{\longrightarrow}
   \Lambda_{2}
   \longrightarrow
   \cdots
\]
For $t\geq 0$, $k\geq 1$ 
we denote by $f_t^{(k)}=f_{t+k} \circ \cdots \circ f_{t+1}$ the $k$'th
iterated map from $\Omega_k(\sigma^t(\FF))$ onto $K$
(with $f_t^{(0)}$ being the identity map on $K$).
We write simply  $f^{(k)}\equiv f^{(k)}_0: \Omega_k(\FF) \rr K$ for the 
iterated map starting at time zero and $Df^{(k)}(x)$ for the conformal
derivative of this iterated map.\\

 For $x\in\Lambda_0$ we write
$x_{j}=f^{(j)}(x)$, $j\geq 0$ for its forward orbit (and similarly
for $u\in\Lambda_0$).
 Using this notation we define the
$n$'th Bowen ball around $x$:
\[
   B_{n}(x) =
 \{ u\in \Lambda_0 : 
       d(x_{j},u_{j}) < \delta_{f_{j+1}}(x_{j}),
             0 \leq j \leq n \}.
\]
and then also the 
 $(n-1,\Delta)$-Bowen ball  around 
$x\in\Lambda_0$:
\[
   B_{n-1,\Delta}(x) =
 \{ u\in B_{n-1}(x) : 
              d(x_n,u_n)<\Delta \}.
\]
Then $f^{(n)}:B_{n-1,\Delta}(x) \rr B(f^{(n)}(x),\Delta)$
is a  uniformly expanding homeomorphism for all $x\in\Lambda_0$,
When $u\in B_{n-1,\Delta}(x)$ we say that $u$ and $x$ 
are $(n-1,\Delta)$-close.
Our hypotheses imply that being $(n-1,\Delta)$-close is  a reflexive
relation
 (not so obvious when $\delta_f(x)$ depends on $x$)
  as is shown in the proof of the following

\begin{Lemma}{\bf [Pairing]}
 For $n\in\NN$, $y,w\in K$ with $d(y,w)\leq \Delta$, the sets 
     $(f^{(n)})^{-1}\{y\}$
   and
     $(f^{(n)})^{-1}\{w\}$
may be paired uniquely into pairs of $(n-1,\Delta)$-close points. 
\label{lemma pairing 2}
\end{Lemma}
Proof: Fix $f=f_n$ and
 let $x\in f^{-1}\{y\}$. By (T1) $f B(x,\delta_{f}(x))$
contains $B(f(x),\Delta)\ni w$. By (T2) there is a unique 
preimage $z\in f^{-1}\{w\}$ at a distance 
$d(z,x)< \delta_{f}(x)\leq \Delta$ to $x$.
We claim that then also $x \in B(z,\delta_{f}(z))$
(which makes the pairing unique and reflexive). If this were not
so then there must be $x_2\in B(z,\delta_{f}(z))\cap f^{-1}\{y\}$
which by (T0) must verify:
 $\delta_f(x)>\delta_f(z)>d(x_2,z)>\delta_f(x_2)$.
Inductively one constructs disjoint sequences
$x_1=x,x_2,\ldots\in f^{-1}\{y\}$, $z_1=z,z_2,\ldots\in f^{-1}\{w\}$
for which $\delta_f(x_1)>\delta_f(z_1)>\delta_f(x_2)>\delta_f(z_2)>\cdots$
and this contradicts finiteness of the degree of $f$.
 Returning to the sequence of mappings we obtain by recursion in $n$
the unique pairing. \Halmos

\begin{Lemma}{\bf [Sub-exponential Distortion]}
\label{lemma T distortion}
There is a sub-exponential sequence,
 $(c_n)_{n\in\NN}$, (depending on the equi-distortion function, $\epsilon$,
 but not on the actual sequence of maps)
 such that for any two points $z,u\neq x$
 which are $(n-1,\Delta)$-close to an $x\in\Lambda_0$
  \[
     \frac{1}{c_n} \leq 
               \frac
               {d(f^{(n)}(u),f^{(n)}(x))} 
	       {d(u,x) \; Df^{(n)}(z)} 
	       \leq c_n
      \mbox{\ \ \ \ and \ \ \ \ }
     \frac{1}{c_n} \leq 
               \frac
               {Df^{(n)}(x)}
	       {Df^{(n)}(z)}
	       \leq c_n
  \]
\end{Lemma}
Proof:
As in Lemma
\ref{lemma distortion}.
 More precisely,  we 
have $\log \,c_n=\epsilon(\Delta)+\epsilon(\Delta/\beta)+
\cdots + \epsilon(\Delta/\beta^{n-1})$. \Halmos\\

For $s\geq 0$, $f\in\EE(\Delta,\beta,\epsilon)$ we define as before 
a transfer operator
    $L_{s,f} : \calM(K) \rr \calM(K)$ by setting:
\begin{equation}
   (L_{s,f} \phi)(y)  \equiv 
      \sum_{x\in f^{-1}\{y\}} 
        \left( {Df(x)} \right)^{-s}
        \phi_x,\ \ \  y\in K,
            \phi \in \calM(K).
\label{transfer op time}
\end{equation}

We write $L^{(n)}_s=L_{s,f_n}\circ \cdots \circ L_{s,f_{1}}$ for the 
$n$'th iterated operator from $\calM(K)$ to $\calM(K)$.
We denote by $\bfone=\chi_{\mbox{}_K}$
 the constant function which equals one on $K$
and as in (\ref{lower upper})  we define
(omiting the dependency on $\FF$ in the notation)~:

  \[
      M_n(s)\equiv \sup_{y\in K} L_s^{(n)} \bfone(y)
      \ \ \mbox{and} \ \ 
      m_n(s)\equiv \inf_{y\in K} L_s^{(n)} \bfone (y) ,
   \]
and then the lower and upper $s$-conformal pressures~:
   \[
   -\infty \leq
  {\underline{P}(s)} \equiv \liminf_n \frac1n \log m_n(s) \ \ \ 
     \leq \ \ \ 
  {\overline{P}(s)} \equiv \limsup_n \frac1n \log M_n(s)
    \leq +\infty.
\]
In general these limits need not be equal nor finite.
Explicitly, we have e.g.\ the formula for the lower pressure,
 similar to (\ref{conf deriv}),
 \[
    \UP(s)=\liminf_n \frac1n \log \inf_{y\in K}
           \sum_{x\in  (f^{(n)})^{-1}\{y\}} 
	    \left( Df^{(n)}_x \right)^{-s}.
 \]

We define the following {\em lower and upper critical exponents}
with values in $[0,+\infty]$:
  \[
       {\underline{s}}_\crit 
           = \sup\{s\geq 0: 
	     {\underline{P}(s)} >0\} \ \ \mbox{and} \ \ 
       {\overline{s}}_\crit
           = \inf\{s\geq 0: 
	     {\overline{P}(s)} <0\}.
  \]

% \begin{Remarks}
% Note that since $\UP(s)+s \; \log \; \beta$ is non-decreasing we have
% that $\UP(s)<0$ for $s>\us_\crit$.
% \end{Remarks}

It will be necessary to make some
additional assumptions on mixing and growth rates.
For our purposes the following suffices:

\begin{Assumption}
\label{T4-T5}
\mbox{}

\begin{itemize}
  \item[(T4)] There is $n_0<\infty$ such that 
      the sequence $(f_k)_{k\in\NN}$ 
   is $(n_0,\Delta)$-mixing, i.e.\
    for any $y\in K$, and $t\geq 0$,
        \[
         (f_t^{(n_0)})^{-1}\{y\} \mbox{\ is $\Delta$-dense in \ }
          (f_t^{(n_0)})^{-1} K.
        \]
  \item[(T5)] The sequence
        $(\lambda(f_k))_{k\in\NNs}$
  is sub-exponential, i.e.
       \[
           \lim_k \frac1k \log \lambda(f_k) = 0.
        \]
%       \ \ $(\|Df_k\|)_{k\in\NNs}$,\ $(\delta_{f_k})_{k\in\NNs}$ \ and
%    $(d(f_k))_{k\in\NNs}$  
%  are all sub-exponential.
\end{itemize}
\end{Assumption}

\begin{Lemma} Assuming (T0)-(T5) we have
(the limits need not be finite):
\begin{eqnarray*}
     \OP(s)&=&\limsup\frac1n m_n(s)=\limsup\frac1n M_n(s)\\
     \UP(s)&=&\liminf\frac1n m_n(s)=\liminf\frac1n M_n(s)
\end{eqnarray*}
\label{lemma lower upper}
\end{Lemma}
Proof:
By a small modification  (notably replacing $\delta_f$ by $\Delta$)
in the last half of the proof of the
Operator bounds - Lemma \ref{Operator bounds} - and making use of
mixing  (T4),
 we deduce similarly to 
(\ref{mM bounds}) that
 \[
       m_{n+n_0}(s) \geq 
         (\|Df_{n+1}\|\cdots \|Df_{n+n_0}\| c_n)^{-s} M_n(s)/2.
 \]
 By Lemma \ref{lemma T distortion} the sequence $c_n$ is sub-exponential.
Due to (T5) and since $n_0$ is fixed,
 $M_n(s)/m_{n+n_0}(s)$ is then of 
sub-exponential growth. Whether finite or not,
the above lim inf's and lim sup's 
agree. \Halmos\\

\begin{Lemma} Assuming (T0)-(T5) we have the following dichotomy: Either
 $\Lambda_0$ is a finite set or $\Lambda_0$ is a perfect set.
\end{Lemma}
Proof: Suppose that
$\Lambda_k$ is a singleton for some $k$ (this happens iff the degrees
of the sequence of mappings is eventually one !). Then also 
$\Lambda_n$ is a singleton for all
$n\geq k$ and
$\Lambda_0$ is a finite set because all the (preceeding)
maps are of finite degree.  
Suppose instead that no $\Lambda_k$ is reduced to a singleton
and let us take $x\in \Lambda(\FF)$ as well as
$n\geq n_0$. Let $z\in \Lambda_n$, $z\neq f^{(n)}(x)$.
By (T1),(T2) and (T4) $z$ must have an $n$'th  pre-image in $\Lambda_0$
distinct from $x$ and
at a distance less than $\beta^{n_0-n}\Delta$ to $x$. Thus, $x$ is 
a point of accumulation of other points in $\Lambda_0$. \Halmos\\

We have the following  (see \cite[Theorem 2.1 and 3.8]{Bar96}
for similar results)~:
\begin{Theorem}
\label{Thm time dep conf rep}
 Let $\Lambda_0$ denote the time-zero conformal repeller for a
sequence of $\EE(\Delta,\beta,\epsilon)$-maps,
$(f_t)_{t\in\NN}$, verifying conditions $(T0)-(T5)$.
Then we have the following inequalities
(note that the first is actually an equality),
 regarding 
dimensions of $\Lambda_0=\Lambda(\FF)$,
\[
       {\underline{s}}_\crit
    = \dimH \Lambda_0 \leq \boxinf \Lambda_0 \leq \boxsup \Lambda_0 \leq 
       {\overline{s}}_\crit.
\]
If, in addition $\lim\frac1n \log m_n(\us_\crit)=0$  
 then $\us_\crit=\os_\crit$ and
all the above  dimensions agree.
\end{Theorem}

Proof: When $\Lambda_0$ is a finite set it is easily seen that
$\OP(0)= 0$ and then that $\us_\crit=\os_\crit=0$ in agreement
with our claim.  In the following
we assume that $\Lambda_0$ has no isolated points.\\

  ($\us_\crit\leq \dimH \Lambda_0$):
Let $U$ be an open subset intersecting $\Lambda_0$
and of diameter not exceeding $\delta_{f_1}$. Choose 
$x=x(U)\in U \cap \Lambda_0$ 
and take again $k=k(U)\geq 0$
to be the largest integer (finite when $\Lambda_0$
is without isolated points)  such that
    $U\subset B_{k}(x)$.
Then there is $u\in U\setminus B_{k+1}(x)\subset B_k(x)\setminus B_{k+1}(x)$
 for which we have
$\delta_{f_{k+2}} \leq d(x_{k+1},u_{k+1})
 \leq \lambda_1(f_{k+1}) d(x_k,u_k)$.
The bound (\ref{Lk bound}) is replaced by
\[
      L^{(k)}_s \chiU \leq (\diam U)^s 
          \left[ \left( \frac{\lambda_1(f_{k+1}) c_k}
	  {\delta_{f_{k+2}}}\right)^s 
           \frac{1}{m_k(s)} \right] L^{(k)}_s \bfone .
\]
By hypothesis (T5) $\lambda_1(f_k)$ is 
a sub-exponential sequence. $\Delta$-homogeneity,
or more precisely the bound (\ref{def k}), shows that
$\delta_{f_k}\geq \kappa / \lambda_1(f_k)$ is also sub-exponential.
If $\us_\crit=0$ there is nothing to show.
If $0\leq s<\underline{s}_\crit$ then $m_k(s)$ tends
 to infinity exponentially fast and
 the factor in the square bracket is uniformly bounded
from above  by a 
constant $\gamma_1(s)<\infty$. Now, let us take the resctriction
of the resulting inequality  to 
the subset $\Lambda_k$. On the left hand side
$U$ is replaced by its intersection with $\Lambda_0$ and on the
right 
 $L^{(k)}_s \bfone$ is replaced 
by $L^{(k)}_s \chi_{\mbox{}_{\Lambda_0}}$ 
(since $(f^{(k)})^{-1} \Lambda_k=\Lambda_0$):
\[
   L^{(k)}_s \chiUL \leq \gamma_1(s) \, (\diam U)^s\,  L^{(k)}_s \chiLz .
\]
We may then repeat the argument from section
\ref{lower bounds}
to  conclude that
$\dimH \Lambda_0 \geq 
       {\underline{s}}_\crit$.\\

 ($\us_\crit\geq \dimH \Lambda_0$):
To obtain this converse inequality,
we will use a standard trick which amounts to construct
explicit covers of small diameter and give bounds for their Hausdorff
measure. 

By our initial assumption $K$ has a finite
$\Delta$-cover $\{V_1,\ldots,V_{N_\Delta}\}$.
Fix $n\geq 1$
 as well as $i\in\{1,\ldots,N_\Delta\}$.
 Pick $x_i\in V_i$ and
write $(f^{(n)})^{-1}\{x_i\}=\bigcup_{\alpha\in I_i} \{x_{i,\alpha}\}$ 
over a  
finite index set $I_i$. By the Pairing Lemma
\ref{lemma pairing 2},
 we see that
to each $x_{i,\alpha}$ corresponds a pre-image 
$V_{i\alpha}=(f^{(n)})^{-1}V_i \cap B_{n-1,\Delta}(x_{i,\alpha})$ 
(the union over $\alpha$ yields a partition of $(f^{(n)})^{-1}V_i$).
Whence,
by sub-exponential distortion, Lemma 
\ref{lemma T distortion},
\[
   \diam\; V_{i,\alpha} \leq  \frac{2 c_n \Delta}
           {Df^{(n)}_{x_{i,\alpha}}}.
\]
Then,
\[
   \sum_\alpha
   (\diam\; V_{i,\alpha})^s \leq (2 c_n \Delta)^s (L_s^n \bfone)(x_i)
\]
and consequently 
\[
   \sum_{i,\alpha}
   (\diam\; V_{i,\alpha})^s \leq [N_\Delta (2 c_n \Delta)^s 
                          M_n(s)].
\]
Let $s>\underline{s}_\crit$.
Then $\underline{P}(s)<0$ and there is a sub-sequence $n_k$, $k\in\NN$
for which $m_{n_k}$ and, by Lemma
\ref{lemma lower upper},
also $M_{n_k}(s)$
 tend exponentially fast to zero.
For that sub-sequence 
the expression in the square-brackets is
 uniformly  bounded in $n_k$.
Since $\diam\; V_{i,\alpha} \leq 2 c_{n} \Delta\; \beta^{-n}$
which tends to
zero with ${n}$ the
family $\{V_{i,\alpha}\}_{n_k}$ exhibits  covers 
of $\Lambda_0$ of arbitrarily
small diameter.
This implies that \HD($\Lambda$) does not exceed $s$,
whence not $\underline{s}_\crit$.\\

($\boxsup \Lambda_0 \leq \os_\crit$): 
For the upper bound on the box dimensions, 
 consider for 
$0<r<\delta(f_1)$, $x\in\Lambda_0$ the ball
 $U=B(x,r)$ and let $k=k(x,r)\geq 1$ be the smallest 
integer such that
$B_{k-1,\Delta}(x)\subset U$.
Then there is $y\in B_{k-2,\Delta}(x)\setminus U\supset 
     B_{k-2,\Delta}(x)\setminus B_{k-1,\Delta}(x)$.  
 As in section
\ref{dimB upper} we deduce that
\[
    L^{(k)}_s\chiU \geq r^s
         (c_{k-1} \Delta\|Df_k\| )^{-s}
	       \chi_{B(f^{(k)}_x,\Delta)} .
\]
Iterating another $n_0$ times we will 
by hypothesis (T4) cover all of $\Lambda_{k+n_0}$. 
Reasoning as in section \ref{dimB upper} it follows that
\[
    L^{(k+n_0)}_s\chiU \geq (4r)^s \left[
         (4 c_{k-1} \Delta \prod_{j=0}^{n_0}
            \|Df_{k+j}\|)^{-s}
	      \frac{1}{M_{k+n_0}(s)} \right] L_s^{({k+n_0})} \chiLz
\]
If $s>\os_\crit$ the sequence $M_k(s)$ tends to zero exponentially
fast. The sub-exponential bounds in hypothesis (T5)
 imply that the factor in the brackets 
 remains uniformly bounded from below.
We may proceed 
to conclude that 
$\boxsup \Lambda$ does not exceed $s$, whence not $\os_\crit$.\\

Finally, for the last assertion suppose that 
$\frac1n\log m_n(\us_\crit)=0$, i.e.\ the limit exists and equals zero
(cf. the Remark below). Then
Lemma \ref{lemma lower upper}
shows that
  the lower
and upper pressure agree and therefore
$\UP(\os_\crit)=\OP(\os_\crit)=0$. Now, both pressure functions are 
strictly decreasing (because $\beta>1$). Therefore, 
$\os_\crit=\us_\crit$  and
the conclusion follows.
\Halmos\\

\begin{Remarks}
A H\"older inequality (for fixed $n$) shows that
 $s\mapsto \frac1n \log M_n$ is convex in $s$. The property
of being convex is preserved 
when taking limsup (but in general not when taking liminf) so that
 $s\mapsto \OP(s)$ is convex.
Even when $\frac1n \log \; M_n(\os_\crit)$ converges, however,
it can happen that  the limit is $+\infty$
for $s<\os_\crit$. In that case convergence of
$\frac1n \log M_n(\os_\crit)$ could be towards a strictly negative number
and $\us_\crit$ could turn out to be strictly smaller than $\os_\crit$.
\end{Remarks}
% \begin{Theorem}
% \label{Thm random}
% Something with randome conf repellers....
% 
% \end{Theorem}

\section{Random conformal maps and parameter-dependency}
The distortion function, $\epsilon$, gives rise to a natural metric on
 $\EE\equiv \EE(\Delta,\beta,\epsilon)$.  We assume in the following
that $\epsilon$ is extended to all of $\RR_+$ and is a strictly increasing
concave function (or else replace it by an extension of its concave `hull'
and make it increasing). For $f,\tf \in \EE$ we set
$\dE(f,\tf)=+\infty$ if there is $y\in K$ for which 
$\# f^{-1} y\neq \tf^{-1}y$.
Note that by pairing $\# f^{-1}y$ is locally constant.
When the local degrees coincide everywhere
we proceed as follows: For $y\in K$,
we let $\Pi_y$ denote the family of bijections,
$\pi: f^{-1} y \rr \tf^{-1} y$, and set
  \[
    \rho_{\pi,x}(f,\tf)
        = \epsilon \left( \frac{\beta}{\beta-1} d(x,\pi(x)) \right)
              + \left|  \log \frac {D\tf \circ \pi(x)}{Df(x)} \right|.
 \]
The distance between $f$ and $\tf$ is then defined as
\begin{equation}
   \dE(f,\tf) = \sup_{y\in K} \ \inf_{\pi\in \Pi_y} \ \sup_{x\in f^{-1}(y)}
             \ \rho_{\pi,x}(f,\tf).
\end{equation}
Our hypotheses on $\epsilon$ imply that
    $\rho_{\pi_2\circ\pi_1,x}(f_1,f_3)\leq 
    \rho_{\pi_1,x}(f_1,f_2) +
    \rho_{\pi_2,\pi_1(x)}(f_2,f_3)$ from which we deduce that
$\dE$ fulfills a triangular inequality. It is then checked that indeed,
$\dE$ defines a metric on $\EE$.

\begin{Lemma}
\label{lemma param}
 Let $u\leq \Delta$ and $\dE(f,\tf)\leq \epsilon(u)$.
Then for all $y,\ty\in K$ with $d(y,\ty)\leq u$ there exists a pairing
$(x_\alpha,\tx_\alpha)$, $\alpha\in J$ (some index set) of
$f^{-1}(y)$ and $\tf^{-1}(\ty)$ for which $\forall \alpha\in J$,
 \[   d(x_\alpha,\tx_\alpha) \leq u  \ \ \ \mbox{and} \ \ \ 
       \left| \log \frac {Df(x_\alpha)}{D\tf(\tx_\alpha)} \right|
            \leq 2\epsilon(u) .\]
\end{Lemma}
Proof: Let $x\in K$ and choose a bijection $\pi:f^{-1}(y)\rr \tf^{-1}(y)$
for which $\epsilon \left( \frac{\beta}{\beta-1} d(x,\pi(x)) \right)\leq 
          \epsilon(u)$ for all $x\in f^{-1}(y)$.
Then $ d(x,\pi(x)) \leq (1-\frac{1}{\beta}) u$, 
           for all $x\in f^{-1}(y)$. Consider $x\in f^{-1}(y)$ and
$x'=\pi(x)$.
As $\tf B(x',\delta_\tf(x')) \supset B(y,\Delta) \ni \ty$ there is
a unique point $\tx\in B(x',\delta_\tf(x'))$ for which $\ty=\tf(\tx)$.
As the association (for fixed $\pi$),
$x \mapsto x'=\pi(x)\mapsto \tx$ is unique we have
obtained a pairing.
By expansion of $\tf$ we have $d(x',\tx)\leq d(y,\ty)/\beta\leq u/\beta$.
Therefore also,
 $d(x,\tx)\leq u(1-\frac{1}{\beta}) + \frac{u}{\beta} = u$ as wanted.
By definition of distortion we have
$|\log D\tf(x')/D\tf(\tx)| \leq \epsilon(d(y,\ty)) \leq \epsilon(u)$.
Also, $\dE(f,\tf)\leq \epsilon(u)$ implies
$|\log Df(x)/D\tf(x')| \leq  \epsilon(u)$  and the last claim follows.
\Halmos\\ 

Given two sequences 
$\FF=(f_n)_{n\in\NN}$ and 
$\tFF=(f_n)_{n\in\NN}$  in $\EE$ we define their  distance
(one could here replace sup by lim-sup),
\begin{equation}
        d_{\infty}(\FF,\tFF)= \sup_n \dE(f_n,\tf_n).
\end{equation}

\begin{Proposition}
When $d_{\infty}(\FF,\tFF)\leq r \leq \epsilon(\Delta)$ then:
 \[  \left| \UP(s,\FF) - \UP(s,\tFF)\right| \leq 2 r s ,
      \ \ \ s\geq 0 \ \ \ \mbox{and}\]
\[  \left( 
       1+ \frac{2r}{\log \beta}
     \right)^{-1} \leq 
     \frac{\us_\crit(\FF)}{\us_\crit(\tFF)}  \leq
       1+ \frac{2r}{\log \beta} .\]
We have the same bounds for the
upper pressures, $\OP$, and upper critical value, $\os_\crit$.
\label{param depend}
\end{Proposition}

Proof: We perform a recursive pairing of pre-images at distances
less than $u$, with $\epsilon(u)\leq r$. By Lemma 
\ref{lemma param}
for the bounds on the derivatives we obtain
  \[
    \frac{1}{k} \left| 
       \log \frac
         {L^{(k)}_{s,\FF} \bfone (y)}
         {L^{(k)}_{s,\tFF} \bfone (y)} \right| \leq 2rs.
 \]
The first claim follows by taking a limit. For the  second claim
suppose that $s_c=\us_\crit(\FF) < \ts_c =\us_\crit(\tFF)$.
Since $s\mapsto \UP(s,\FF)+s\log \beta$ is non-increasing
(same for $\tFF$) we have
$(\ts_c-s_c)\log \beta \leq \UP(s,\tFF)-\UP(s,\FF) \leq 2 rs$
for all $\ts_c\geq s \geq s_c$. From this inequality the 
other bound follows.\Halmos\\

We associate to the metric space $(\EE,dE)$ its corresponding
Borel $\sigma$-algebra and this allows us to construct measurable
maps into $\EE$.  
In the following,
 let $(\Omega,\mu)$ be a probability space and $\tau:\Omega\rr\Omega$
a $\mu$-ergodic transformation.

\begin{Definition} We write 
    $\EE_\Omega\equiv \EE_\Omega(\Delta,\beta,\epsilon)$ for the 
  space of measurable maps,
      $\bff:\omega\in (\Omega,\mu) \mapsto \bff_\omega \in (\EE,\dE)$,
   whose image is almost surely
separable (i.e.\ the image of a subset of full measure
 contains a countable dense set).
Following standard conventions
we say that the map is {\em Bochner}-measurable.
\label{measurable family}
\end{Definition}

We write 
$\FF_\omega=(\bff_{\tau^{n-1}\omega})_{n\in\NN}$ 
for the sequence
of maps fibered at the orbit of $\omega\in\Omega$.
 Denote by $\bff^{(n)}_\omega= \bff_{\tau^{n-1}(\omega) } \circ
 \cdots \circ \bff_\omega $ the iterated map
defined on the domain 
$\Omega_n(\FF_\omega)= \bff_\omega^{-1} \circ \bff_{\tau(\omega)}^{-1} \circ
 \cdots \circ \bff_{\tau^{n-1}(\omega)}^{-1}(K)$.
 The `random' Julia set, as before,
is the compact, non-empty  intersection
     \begin{equation}
       J(\FF_\omega) = \bigcap_{n\geq 0} \Omega_n(\FF_\omega).
        \label{random Julia set}
      \end{equation}

Our assumptions imply that
 $(f_1,\ldots,f_n)\in \EE^n \mapsto f_1^{-1}\circ \cdots \circ
 f_n^{-1}(K) \subset K$  is continuous
 when $K$ is equipped with the Hausdorff topology
 for its non-empty subsets.
It follows that
$\omega \mapsto \Omega_n(\FF_\omega)$ is measurable.
 Uniform contraction implies that $\Omega_n$ convergences exponentially
fast to $J(\FF_\omega)$ in the Hausdorff topology,
whence the `random' conformal repeller, $J(\FF_\omega)$ is (a.s.) 
measurable for the Hausdorff $\sigma$-algebra.

Using the estimates from the previous Proposition,
the function,
 $(f_1,\ldots,f_n) \in \EE^n \mapsto M_n(s,(f_1,\ldots,f_n))$
is continuous. Almost sure separability of 
$\{\bff_\omega:\omega\in\Omega\}\subset \EE$ implies then that
$\omega\mapsto M_n(s,\FF_\omega)$ is measurable (with the standard
Borel $\sigma$-algebra on the reals). For example, if
$V_1,V_2$ are open subsets of $\EE$, the pre-image of $V_1\times V_2$
by
$\omega\mapsto (\bff_\omega, \bff_{\tau\omega})$ is
$\bff^{-1}(V_1) \cap \tau^{-1}\bff^{-1}(V_2)$ which is
measurable.
The function $\OP(s,\FF_\omega)$, being a lim sup of measurable
functions, is then also measurable
 (and the same is true for $m_n$ and $\UP$).
We  define the distance between 
 $\bff, \tbff \in \EE_\Omega$ to be
\begin{equation}
    d_{\EE,\Omega} (\bff,\tbff) =
       \mu\mbox{-ess} \, \sup_\omega 
      d_{\EE}(\bff_\omega,\tbff_\omega) \in [0,+\infty].
\end{equation}

\begin{Theorem}
\label{Thm time depend}
Let $\tau$ be an ergodic transformation on $(\Omega,\mu)$ and
let $\bff=(f_{\omega})_{\omega\in\Omega}\in \EE_\Omega$ be 
Bochner-measurable (Definition \ref{measurable family}).
We suppose that there is $n_0<\infty$ such that
 almost surely the sequence
$\FF_\omega=(\bff_{\tau^{n-1}\omega})_{n\in\NN}$ 
 is $(n_0,\Delta)$-mixing
(Condition (T4) in Assumption \ref{T4-T5}).\\

(a) Suppose that $\EEE\log \|D\bff_\omega\| < + \infty$.  
[We say that the family is of {\em bounded average logarithmic dilation}].
Then for any $s\geq 0$ and
$\mu$-almost surely, the pressure function
$\UP(s,\FF_\omega)$ is independent of $\omega$. We write
$\UP(s,\bff)$ for this almost sure value.
The various dimensions of the random conformal repeller agree
(a.s.) in value. 
Their common value  is  (a.s.) constant and given by
\[
  \dim \Lambda(\FF_\omega) =
  \sup \{s\geq 0: \UP(s,\bff)>0 \} \in [0,+\infty].
\]
\\

(b) The (a.s.) dimension is finite iff \, $\UP(0,\bff)< +\infty$
(this is the case, e.g. if
  \,$\EEE \log d^o_{\max}(\bff) < \infty$) 
 and  one has the estimate,
\[ 
 \frac{\EEE \log d^o_{\min}(\bff)} {\EEE \log \|D\bff\|}
  \leq \dim \Lambda(\FF_\omega) \leq
 \frac{\EEE \log d^o_{\max}(\bff)} {-\EEE \log \|1/D\bff\|}
 .\]

(c) The mapping,
   $\bff \in (\EE_\Omega,d_{\EE,\Omega}) \mapsto \dim \Lambda(\FF_\omega)$,
   is $\frac{2}{\log \beta}$-Lipschitz 
(at distances  $\leq \epsilon(\Delta)$).
\end{Theorem}

Proof:

Write $\phi=\phi_\omega = \log \| D\bff_\omega \| \geq 0$ and similarly
$\phi^{(n)}=\phi^{(n)}_\omega=\log \| D\bff^{(n)}_\omega\| \geq 0$.
Then $\phi^{(n)}\leq \phi^{(k)} + \phi^{(n-k)}\circ \tau^k$, $0<k<n$
       and since $\phi$ is integrable we get by Kingman's subergodic
Theorem, \cite{King68}, that the limit
  \[
     \lim_n \frac1n \phi^{(n)} \geq 0
  \]
exists $\mu$-almost surely. As a consequence,
  \[
     \lim_n \frac1n \phi\circ\tau^n =
     \lim \frac{n+1}{n} \; \frac{1}{n+1} \phi^{(n+1)}
             - \frac{1}{n} \phi^{(n)} = 0
  \]
$\mu$-almost surely. Thus the sequence of maximal dilations is 
almost surely sub-exponential  (Condition (T5) of 
Assumption \ref{T4-T5}). Condition (T4) of that assumption is
a.s.\ verified by the hypotheses stated in our Theorem.
It follows by Theorem  
\ref{Thm time dep conf rep} that the Hausdorff dimension of the
random repeller, $\Lambda(\FF_\omega)$ a.s.\ is given by
 $\us_\crit(\FF_\omega)$.
We wish to show that a.s.\ the value is constant and
that a.s.\ $\frac1n \log m_n(\us_\crit(\FF_\omega)) \rr 0$
as $n\rr\infty$.

We have the following bounds for the action of the transfer operator, $L_{s,f}$,
upon a positive function, $\phi> 0$~:
\begin{equation}
\frac{d^o_{\min}(f)}{\|Df\|^s} \; \min \phi 
\leq L_{s,f} \phi \leq
{d^o_{\max}(f)}\; {\|\frac{1}{Df}\|^s}\;  \max \phi .
\label{op bounds}
\end{equation}

Here, $d^o_{\max}(f)$ and $d^o_{\min}(f)$ denotes the maximal,
respectively, the minimal (local) degree of the mapping $f$.
From the lower bound we obtain in particular,
\[ 
\EEE \log m_1(s,\FF_\omega) \geq
   \EEE \log d^o_{\min}(\bff) - s \EEE \log \| D \bff\| \geq
    - s \EEE \log \| D \bff\|.\]

 The family, $m_n$,    is super-multiplicative, i.e.\ 
 $m_n(s,\FF_\omega) \geq m_{n-k}(s,\FF_{\tau^k\omega}) m_k(s,\FF_\omega)$,
for $n,k\geq 0$ and $\omega\in\Omega$.
Writing $\log_x x = \max \{ 0, \log x \}$, $x>0$,
we have

\[ \EEE \log_+ \frac{1}{m_1(s,\FF_\omega)} \leq 
     s \; \EEE \log \| D\bff \| .
\]
As the latter quantity is assumed finite, we may apply Kingman's
super-ergodic Theorem to $m_n$ (i.e.\ the 
sub-ergodic Theorem to the sequence $1/m_n$), to deduce that the limit
\[ \frac1n \lim \log m_n(s,\FF_\omega) \in (-\infty,+\infty]\]
exists $\mu$-almost surely and is a.s.\ constant. We write $\UP(s,\bff)$ for
this a.s.\ limit. From the expression for the operator 
and for fixed $n$ and $\omega\in\Omega$, the sequence,
$\|D\bff_\omega^{(n)}\| m_n(s,\FF_\omega)$,
is a non-decreasing function of $s$. The same is then true for
\[ \frac1n \log \| D\bff^{(n)}_{\omega}\| +
   \frac1n \log m_n(s,\FF_\omega) .\]
Apply now Kingman's sub-ergodic, respectively super-ergodic, 
 Theorem to these two terms to see that
\[ s \, \EEE\, \log \|D\bff\| + \UP(s,\bff) \in (-\infty,+\infty]\]
is a non-decreasing function of $s$. It is seen in a similar way that
\[ s\, \log \beta + \UP(s,\bff)\in (-\infty,+\infty]\]
is  non-increasing.
These two bounds
 together 
with Theorem
\ref{Thm time dep conf rep}
 imply that either 
(1) $\UP(0,\bff)=\infty$,
 $\UP(s,\bff)$ is infinite for all $s\geq 0$ and
$\us_\crit=+\infty$, or
(2) $\UP(0,\bff)< +\infty$ in which case the function
 $s\mapsto \UP(s,\bff)$ is continuous, strictly decreasing and
has a unique zero $\us_\crit$.
The additional condition in Theorem
\ref{Thm time dep conf rep} is thus satisfied
 and $\us_\crit$ therefore equals all of the various dimensions.
The estimate, (b), for the dimensions follows from
(\ref{op bounds}) and taking averages as above.
Finally, (c) is a consequence of 
Proposition \ref{param depend} and the fact that 
$\us_\crit$ a.s.\ equals the dimensions.
\Halmos\\

\begin{Example}
Let $K=\{\phi\in \ell^2(\NN): \|\phi\|\leq 1\}$ and denote
by $e_n$, $n\in \NN$ the canonical basis for $\ell^2(\NN)$.
The domains $D_n=\Cl \; {B(\frac23 e_n,\frac16)}$,
$n\in \NN$ maps conformally onto $K$ by
$x\mapsto 6(x-\frac23 e_n)$.
 To each $n\in\NN$ we consider the conformal map,
$f_n$, of degree $n$, which maps $D_1\cup \ldots\cup D_n$ onto $K$
by the above mappings. Finally let $\nu$ be a probability measure
on $\NN$. Picking an i.i.d.\ sequence of the mappings, $f_n$,
 according to
the distribution $\nu$ we obtain a conformal repeller for which all
dimensions almost surely agree.
In this case we have equality in the estimates in
Theorem \ref{Thm time depend} (b) so
 the a.s.\ common value for the dimensions is given by
\[
   \frac{\sum_n n \; \nu(n)}{\log 6}.
\]
Finiteness of the dimension thus depends on $n$ having finite 
average or not, cf.\ also \cite[Example 2.1]{DT01}.

\end{Example}

The Lipschitz continuity of the dimensions with respect to
parameters is somewhat delusive because it is with respect
to our particular metric on $\EE$. In practice,
when constructing parametrized families of mappings it is really
the modulus of continuity of $Df$, i.e.\ the $\epsilon$-function
in $\EE(K,\Delta,\epsilon)$ that comes into play~:

\begin{Example}
We consider here just the case of one stationary map, $f\in\EE$.
Let $T_t$, $t\geq 0$
 be a Lipschitz motion of $\Omega_f$. By this we mean that
$T_t^{-1}:\Omega_f \rr K$, $t\geq 0$, is a family
of conformal mappings  with $T_0=\mbox{id}$, 
$|\log DT_t^{-1}(x)|\leq t$, $x\in\Omega_f$ and
 $\sup_{x\in K} d(x,T_t^{-1}x) \leq t$. Let
$\epsilon_{T_t}(r)$ denote the distortion function for $T_t$
(which we may defined in the same way as for $f$ when $r<\Delta-t$). 
 A calculation then shows that
for $t$ small enough, $\epsilon_{f\circ T_t}(r) \leq
       \epsilon_f(r) + \epsilon_{T_t} (r/\beta)$. One 
also checks that
$d_\EE(f\circ T_t,f)
 \leq 2\epsilon_f(t)+ \epsilon_{T_t}(t/\beta)+ c t$.
By Theorem \ref{Thm time depend} (c),
 the mapping  $t\mapsto d(t)= \dimH \Lambda(f\circ T_t)$
for $t$ small verifies
     \[ |\log \frac{d(t)}{d(0)} | \leq
        \frac{2}{\log \beta} ( 2 \epsilon_f(t) + 
                   \epsilon_{T_t}(\frac{t}{\beta}) + c t).
    \]
When Thermodynamic Formalism applies, in particular when a
bit more smoothness is imposed, a similar result could be deduced
within the framework (and restrictions) of TF. 
I am, however, not aware of any 
results published on this.
\end{Example}

\newpage

\section{Part II: Random Julia sets and parameter dependency}
 Let
$U\subset \widehat{\CC}$ be an open non-empty  connected 
subset of the Riemann sphere 
 omitting at least three points. 
 We denote by
$(U,d_U)$ the set $U$ equipped with a hyperbolic metric.
 As normalisation we use
 $ds = 2|dz|/(1-|z|^2)$ on the unit disk $\DD$ and the hereby induced metric
for the hyperbolic Riemann surface $U$
(cf.\ Remark \ref{Riemann surfaces} below).
 In particular, we have
for the unit disk and $z\in\DD$,
  \[
         d_{\DD}(0,z) =  \log \frac{1+|z|}{1-|z|}, \ \ \ \
        |z| = \tanh \frac{d_\DD(0,z)}{2}.
  \]
We write $B(u,r)\equiv B_U(u,r)$ for the hyperbolic ball
 of radius $r>0$ centered
at $u\in (U,d_U)$, $B_\DD(t,r)$ for the similar hyperbolic 
ball in $(\DD,d_\DD)$
 and $B_\CCs(u,r)=\{z\in\CC: |z-u|<r\}$
 for a standard Euclidean ball in $\CC$. 

Recall that when $K\subset U$ is a compact subset 
the inclusion mapping $(\Int K,d_{\Int K}) \hookrightarrow (\Int K,d_U)$
is a strict contraction  
\cite[Theorem 4.2]{CG93}
by some factor $\beta=\beta(K,U)>1$, depending on $K$ and $U$ only.
  We consider
the family $\EE(K,U)$ of finite degree unramified conformal covering maps
\[
f:\D_f\rr U
\]
for which 
the domain $\D_f$ is a subset of the compact set $K$.
We may assume without loss of generality that $K$ is the
closure of its own interior.
Our first goal is to show that such maps {\em a fortiori}
verify conditions (T0)-(T3) from the previous section, in which 
the set $K$ is the same as here
 and the metric $d$ on $K$ is the restriction of 
the hyperbolic metric $d_U$ to $K$.\\

\begin{figure}
\begin{center}
\epsfig{figure=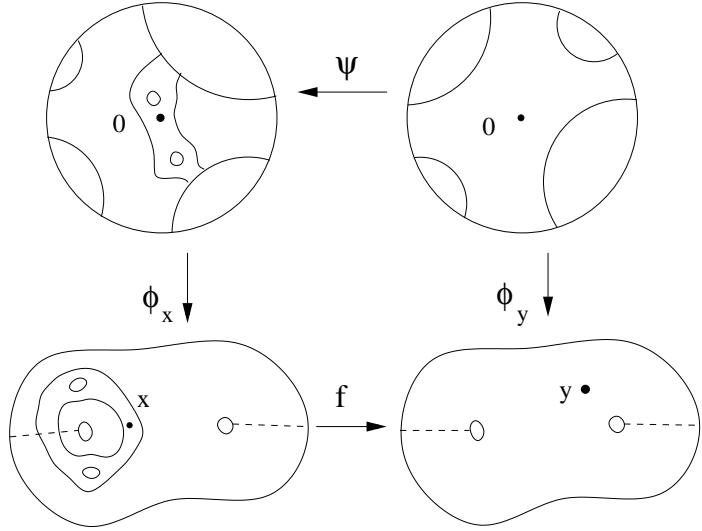,height=7cm}
\label{ramcov}
\end{center}
\caption{An example of
     a covering map of degree 2 and its `inverse' in the universal cover.
       Cuts along the dotted lines become arcs in the lift.
  One fundamental domain is sketched in each cover.}
\end{figure}

Let $\ell=\ell(K,U)>0$ be the infimum length of closed non-contractible
curves (sometimes called essential loops)
 intersecting $K$ and
let $\alpha= \tanh (\ell/4)$ 
( $\ell=+\infty$ and
$\alpha=1$ 
when $U$ is simply connected).
We define constants
\begin{equation}
     \Delta=\Delta(K,U) = \log \frac{1 + \alpha/7}{1-\alpha/7} , \ \ \ 
     \Delta'=\Delta'(K,U) = \log \frac{1 + \alpha/2}{1-\alpha/2} 
     \label{def Delta}
\end{equation}

and for $0\leq r < \ell/2$ the $\epsilon$-function 
\begin{equation}
       \epsilon_\ell(r)=
        - 6 \log 
                \left(  1-\frac{\tanh(r/2)}{\tanh(\ell/4)}   \right).
       \label{epsilon function}
\end{equation}
One has~: 
$\tanh \frac{\Delta}{2} = \frac{\alpha}{7}$,
$\tanh \frac{\Delta'}{2} = \frac{\alpha}{2}$,
$\Delta<\ell/14$,
$\Delta'<\ell/4$ and 
$\epsilon_\ell(\Delta)=6 \log 7/6 <1$.

\begin{Remarks}
\label{Riemann surfaces}
 We recall some facts about universal covering maps
of Riemann surfaces:
         Let $\phi:\DD\rr U$ be a universal conformal covering map of $U$.
For $x,y\in U$ their hyperbolic distance are defined as
 $d_U(x,y) =\min \{d_\DD(\hatx,\haty)\}$
where the minimum is taken over lifts
 $\hatx \in \phi^{-1}(x)$ and $\haty \in \phi^{-1}(y)$  of $x$ and $y$,
respectively.
If $p$, $p'\in \phi^{-1}(y)$ are two disctinct lifts
 of a point $y\in K$ then $d_\DD(p,p')\geq \ell$.
Otherwise the geodesic connecting $p$ and $p'$ projects to a closed
non-contractible  curve
in $U$ intersecting $K$ and of length $< \ell$, contradicting 
our definition of $\ell$. For the same reason,
 the map $\phi:B_\DD(p,\ell/2) \rr B(y,\ell/2)$ 
must be a conformal bijection 
which preserves distances to $y$,
i.e.\ if $z\in B_\DD(p,\ell/2)$ then $d_\DD(z,p)=d_U(\phi(z),y)$
Note, however, that $\phi$ need not be an isometry on the full disc, 
since two points in 
$B(y,\ell/2)\setminus K$
 may have lifts closer than their lifts in $B_\DD(p,\ell/2)$.\\
\end{Remarks}

We have the following 
\begin{Lemma}[Local Koebe Distortion]
\label{Local Koebe Distortion}
 Let $f\in\EE(K,U)$. Denote by $\|Df\|$ the maximal 
 conformal derivative of $f$ on the set $f^{-1}K$.
We define $\lambda_1(f)=3 \| Df\|$.
Let $x\in \D_f\cap f^{-1}K$ 
 and set 
 \begin{equation}
    \delta_f(x) = 
          \min\{     \log \frac{5+\alpha/ Df(x)}
              {5-\alpha/ Df(x)} , \Delta \}.
     \label{def delta}
\end{equation}
Let also
  $\delta_f= \min \{\log \frac{5 +\alpha/\|Df\|}
     {5-\alpha/\|Df\|}, \Delta\}>0$
be the minimum value of $\delta_f(x)$
over the compact set $\D_f\cap f^{-1}K$. 
Then  $B(x,\delta_f(x))\subset \D_f$ and
we have the following properties~:
 \begin{enumerate}
   \item[(0)] \ \ If $x'\neq x$ is another pre-image of $f(x)$ then
           $B(x,\delta_f(x))$ and $B(x',\delta_f(x'))$ are disjoint.
   \item[(1)] \ \
     $f$ is univalent on the hyperbolic disk $B(x,\delta_f(x))$
      and $B(f_x,\Delta)\subset fB(x,\delta_f(x))$.
   \item[(2)] \ \ $\beta\leq f[u,x]\leq \lambda_1(f)$ 
        for  $u\in B(x,\delta_f(x))$.
   %\item[(2)] \ \ $\beta\leq f[u,v]\leq \lambda_1(f)$ 
        %for $u,v\in K$, $d(u,v)\leq \min\{\delta_f(u),\delta_f(v)\}$.
   \item[(3)] \ \ If $u,v\in B(x,\delta_f(x))$ and $f_u,f_v\in B(f_x,r)$
     with $0<r \leq\Delta$ then
     \begin{equation}
       \left| \log \left( \frac{d(f_x,f_u)}{d(x,u) Df(v)} \right) \right|
           \leq \epsilon_\ell(r).
            \label{epsilon bound}
    \end{equation}
   \end{enumerate}
\end{Lemma}

Proof:  
 Let $C$ be a connected component
of $\D_f\subset K$ and fix an $x\in C$ for which $y=f(x)\in K\subset U$.
Pick universal conformal covering maps,
$\phi_x : \DD \rr U$ and $\phi_y: \DD \rr U$ for which 
$\phi_x(0)=x$ and $\phi_y(0)=y$.
Let $\widehat{C}=\phi_x^{-1} C \subset \DD$  be the lift of the
connected component $C$ containing $x$. The composed map, 
$f\circ \phi_x: \widehat{C} \rr U$ is a conformal covering map of $U$.
Since  $\phi_y:\DD\rr U$ is a universal covering there
is a unique, {\em a fortiori} conformal,
 map $\psi=\psi_{x,y}:\DD \rr \widehat{C}$ such that
$\psi_{y,x}(0)=0$ and (cf. figure \ref{ramcov}),
\[
   f \circ \phi_x \circ \psi_{x,y} \equiv \phi_y : \DD \rr U .
\]
By definition of the hyperbolic metric
the conformal derivative of $f$ at $x$ is given by
\[
   \lambda \equiv  Df(x) = 1/|\psi'(0)|.
\]
More generally, if $u=\phi_x\circ \psi(z) \in C$, $z\in\DD$ then 
 \[
    \label{conform deriv}
    Df(u)= 1/D\psi(z) =\frac1{|\psi'(z)|} \frac{1-|\psi(z)|^2}{1-|z|^2}.
\]
The value does not depend on the choices of covering maps
 because the conformal line element $ds=2|dz|/(1-|z|^2)$ is invariant under
conformal automorphisms of the unit disk (both in the source and in
the image).

The map $\psi: (\DD,d_\DD) \rr (\widehat{C},d_{\widehat{C}})$ is non-expanding
\cite[Theorem 4.1]{CG93}. As mentioned above the inclusion 
$(\widehat{C},d_{\widehat{C}}) \hookrightarrow (\widehat{C},d_U)$ is 
$\beta^{-1}$-Lipschitz so the composed map 
$\psi: (\DD,d_\DD)\rr (\DD,d_\DD)$ is also $\beta^{-1}$-Lipschitz.

The map $\psi$ need not, however, be univalent on all of $\DD$, because
a non-contractible loop in $C$ may be contractible in $U$
(as is the case in figure \ref{ramcov}).
On the other hand, the map 
$\phi_y:B_\DD(0,\ell/2) \rr B(y,\ell/2)$ is a conformal
bijection
(Remark \ref{Riemann surfaces})
 so that 
 \[
      h = \phi_x \circ \psi \circ \phi_y^{-1} : B(y,\ell/2)
           \rr B(x,\ell/ (2 \beta))
\]
defines a local inverse of $f$.
 In particular, we see that $\psi$ is univalent on
the disk $B_\DD(0,\ell/2)=B_\CCs(0,\alpha)$.
 The map, $g:\DD\rr\CC$, given by
     \[
      g(t) = \frac{\psi(t \alpha)}{\alpha \psi'(0)},
     \]
is therefore univalent and normalised so that
$g(0)=0$ and $g'(0)=1$.
 The Koebe distortion Theorem
\cite[Theorem 1.6]{CG93} applied to $g$ 
shows that if $|z| <\alpha$ then
\newcommand{\ds}{\displaystyle}
\begin{equation}
    \begin{array}{ccccc}
         \ds \frac{1}{(1+|z|/\alpha)^2}               &\leq&
         \ds \left| \frac{\lambda \psi(z)}{z} \right|    &\leq&
         \ds \frac{1}{(1-|z|/\alpha)^2} ,   \\[4mm]
         \ds \frac{1-|z|/\alpha}{(1+|z|/\alpha)^3}               &\leq&
         \ds \left| {\lambda \psi'(z)}  \right|   &\leq&
         \ds \frac{1+|z|/\alpha}{(1-|z|/\alpha)^3}    .
     \end{array}
     \label{Koebe bounds}
\end{equation}
Using the first bound one verifies that
\begin{equation}
   \psi B_\CCs(0,\frac{\alpha}{7})
       \subset B_\CCs(0,\frac{ \alpha}{5  \lambda})
         \subset \psi B_\CCs(0,\frac{\alpha}{2}) .
   \label{inclusion}
\end{equation}
Going back to hyperbolic distances and $U$, and noting that also
$\psi B_\CCs(0,\alpha/7)\subset B_\CCs(0,\alpha/7)$,
we obtain 
\[
   h B(y,\Delta) \subset B(x,\delta_f(x)) \subset h B(y,\Delta')
          \subset h B(y,\ell/4) .
\]
     with the definition of $\Delta, \Delta'$
      and $\delta_f(x)$ as in (\ref{def Delta}) and (\ref{def delta}). 
In particular, $B(x,\delta_f(x))\subset C \subset \D_f$.

Property (0):
Let $x'$ be another pre-image of $y$ distinct from $x$.
Since $B(x,\delta_f(x))\subset C$ the balls $B_1=B(x,\delta_f(x))$ and
$B_2=B(x',\delta_f(x'))$ are disjoint if they are in different 
connected components of $\D_f$. If $B_1\cap B_2$ is non-empty 
then we may find a shortest path,
$\gamma \subset B_1\cup B_2$ connecting $x$ and $x'$ within $C$. 
Then $f(\gamma)$ is
a closed non-contractible curve in $U$, containing $y$ and of length 
$< \ell/2+\ell/2=\ell$ which contradicts the definition of $\ell$. 

Property (1): Set $B=B(x,\delta_f(x))$.
    The first inclusion in (\ref{inclusion}) shows that
   $fB \supset B(f(x),\Delta)$ and since the local inverse
   $h$ is well-defined and its image contains $B$ the map
   $f$ is univalent on $B$.

Property (2), $f[u,x]\geq \beta$:\ \ 
 For $v\in B(y,\ell/2)$ we
   have that $d(h(v),h(y))\leq \beta^{-1} d(y,v)$ and since
   $h B(y,\ell/2) \supset B(x,\delta_f(x))$ we see that
   $u=h(v)\in f^{-1}\{v\}$ is the point closest to $x \in f^{-1}\{f(x)\}$.
    Therefore $d(u,x)\leq \beta^{-1} d(y,v)$ and
    we obtain the wanted inequality.

Property (2), $f[u,x]\leq  \lambda_1(f)$:\ \ 
 By Schwarz' Lemma, 
$|\psi(z)|\leq |z|$, $z\in \DD$ and from the expression for the hyperbolic
metric,
 \[
           1 \leq 
          \frac{d_\DDs(z,0)}{d_{\DDs}(\psi(z),0)}
          \frac{|\psi(z)|}{|z|}
            \leq \frac{1}{1-|z|^2},
	   \ \ \ z\in \DD.
  \]
Using the first bound in
 (\ref{Koebe bounds}) we get for $|z|<\alpha$,
\begin{equation}
     {(1- |z|/\alpha)^2}
     \leq  \frac{f[u,x]}{Df(x)} =
    \frac{d_\DDs(\psi(z),0) \lambda}
        {d_\DDs(z,0)} 
        \leq 
     \frac{(1+ |z|/\alpha)^2}{1-|z|^2}.
     \label{fux one}
\end{equation}
In particular, 
for $|z|\leq \alpha/2$ (corresponding to the hyperbolic radius 
$\Delta'$),
\[
   f[x,u] \leq  \frac{(3/2)^2}{1-1/4} \; Df(x) = 3\; Df(x)
    \leq \lambda_1(f),\ \ \
          x\in B(x,\delta_f(x)).
\]
\mbox{}\\[2mm]
\noindent Property (3):\ \ The second bound in 
(\ref{Koebe bounds})  shows that for
 $|z|,|u|\leq \hatr<\alpha$,
  \[ \frac
       {1-\hatr/\alpha}
       {(1+\hatr/\alpha)^3} (1-r^2) \leq
         \frac{Df(x)}{Df(v)} = \lambda |\psi'(u)|
             \frac{1-|u|^2}{1-|\psi(u)|^2} \leq
      \frac {1+\hatr/\alpha}
       {(1-\hatr/\alpha)^3}.\]
Multiplying this and the inequality (\ref{fux one}) we obtain
  \[
    \left| \log \frac {f[u,x]}{Df(v)} \right|
            \leq \log 
          \frac{(1+\hatr/\alpha)^3}{
		  (1-\hatr/\alpha)^3
               (1-\hatr^2)^2}
   \leq 6 \log \frac{1}{1-\hatr/\alpha},
  \]
i.e.\ (\ref{epsilon bound}) with the $\epsilon_\Delta$ function defined
in (\ref{epsilon function}). \Halmos\\

This hyperbolic Koebe Lemma implies that 
conditions (T0)-(T3)
of the previous section are verified for our
class of maps, $\EE(K,U)$, when setting $\Omega_f=\D_f\cap f^{-1}K$ 
and looking at 
the metric space $(K,d_U)$, the $\epsilon$-function 
$\epsilon_\ell$,  and finally $\beta$, $\Delta$ and $\delta_f(x)$
as defined above.\\

\begin{Theorem}
\label{Thm erg analytic}
Let $\tau$ be an ergodic transformation on $(\Omega,\mu)$.
Let $\FF=(f_{\omega})_{\omega\in\Omega}\in \EE_\Omega(K,U)$ be a measurable
family satisfying $\EEE(\log \|Df_\omega\|) < + \infty$.
Then $\mu$-almost surely the various dimensions agree and is given
as the unique zero of the pressure function $P(s)$.
\end{Theorem}
Proof: 
We will  apply Theorem
\ref{Thm time depend}. The assumption of
{\em bounded average logarithmic dilation} is included in our
hypothesis. We need to show that $(n_0,\Delta)$ mixing holds
for some $n_0$.
This follows, however, directly from connectivity of 
$U$ and the properties of our conformal maps.
The diameter of $K$ is finite
within $U$. Given two points $y$ and $z$ in $U$ choose a path of uniformly
bounded length (say less than twice the diameter of $K$) connecting them.
By taking preimages we obtain 
paths of exponentially
shrinking lengths.
It suffices to take $n_0$ such that $2 \,\diam\, K/\beta^{n_0}
 \leq \Delta$
 and (T4) of Assumption \ref{T4-T5} follows.
 An area estimate yields 
$d(f) \Area(K) = \int_{f^{-1}K}|Df|^2 d\Area \leq \|Df\|^2 \Area(K)$, whence
  \begin{equation}
         d(f) \leq \|Df\|^2.
         \label{bd degree}
  \end{equation}
Therefore, $\log d(f)$ is bounded on average and we may apply 
Theorem 
\ref{Thm time depend} to obtain the desired conclusion.
For $\phi\geq 0$ we also have by change of variables,
\[
  \int_K L_{s=2} \phi \; dA =
  \int_{f^{-1}K} \phi \; dA \leq \int_K \phi\; dA,
  \]
which incidently shows that $\os_\crit\leq 2$ (as it should be !).
\Halmos\\

\section{Mirror embedding and real-analyticity of the Hausdorff dimension}
\label{section real analytic}
The dependence of the Hausdorff dimension on parameters may be
studied through the dependence of the pressure function on those parameters.
A complication arise, namely that
 our transfer operators do not
depend analytically on the expanding map. In \cite{Rue82}, Ruelle
circumvented this problem in the case of a (non-random)
 hyperbolic Julia set 
by instead looking at an associated dynamical zeta-function.
 Here, we shall introduce a 
{\em mirror embedding}
which tackles the problem 
directly. We embed our function space into a 
larger space and semi-conjugate
our family of transfer operators  to 
operators with an explicit
real-analytic dependency on parameters and mappings.
We establish a Perron-Frobenius 
theorem  through the contraction
of cones of `real-analytic' functions.
The pressure function may then be calculated as the
averaged action of the operator on a hyperbolic fixed point 
(cf.\ \cite{Rue79,Rue97})
which has the wanted dependence
on parameters. Finally as the pressure function cuts the horisontal axis
transversally the result will follow from another implicit function
Theorem.

\subsection{Mirror extension and mirror embedding}
Let $U$ be a hyperbolic subset of $\HCC$ as before.
We write $\barU=\{\barz: z\in U\}$ for the complex conjugated domain
(not the closure) and we define the mirror extension of $U$ as
the product $\hatU=U \times \barU$.  The map
$j: U \rr \hatU$ given by  $j(z)=(z,\barz)$, $z\in U$ is
a smooth embedding of $U$ onto the mirror diagonal,
   \[
       \diag (U) = \{ (z,\barz): z\in U\}.
    \]
The `exchange-conjugation',
 \[
   c(u,v)=(\overline{v},\overline{u}), \ \ \ (u,v)\in \UbarU
\]
defines an involution on the mirror extension 
leaving invariant the mirror diagonal.
Let $X\subset \hatU$ be an open subset.
We call $X$ mirror symmetric, if
$c(X)=X$.
We say that $X$ is connected to the diagonal if
any connected component of $X$ has a non-empty intersection with $\diag U$.
We write $ A(X)=C^0(\Cl\; X)\cap C^\omega(X)$ for the space 
of holomorphic functions on the mirror extension
 having a continuous extension to the boundary.

\begin{Lemma} Let $X\in \hatU$ be an open, mirror symmetric subset,
connected to the diagonal and let $A=A(X)$. Then 

\begin{enumerate}
\item[(1)] $A$ is a unital Banach algebra
(in the sup-norm) with a complex involution,
\[
   \phi^*(u,v)=\overline{\phi(\overline{v},\overline{u})}
      \equiv \overline{\phi}(v,u),\ \ \ (u,v)\in X, \ \phi\in A.
\]
\item[(2)]
Denote by $A_{\RRs}=\{\phi\in A : \phi^*=\phi\}$,
the space of self-adjoint elements in $A$. Such functions are 
 real-valued on the mirror diagonal and we have
$A=A_{\RRs} \oplus i A_{\RRs}$.
\item[(3)]
       A function $\phi\in A$ is uniquely determined by its restriction
to $(\diag U) \cap X$.\\
\end{enumerate}
\label{unique KD}
\end{Lemma}

Proof: (1) and (2) are clear.
 Suppose now that $\phi$ vanishes on the mirror diagonal.
Because any point in $X$ is path-connected to the diagonal
it suffices to show that $\phi$ vanishes on an open neighborhood
of a diagonal point $(y,\bary)$, $y\in \Int K$. 
For $u,v$ small enough we have a convergent power series expansion
\[
   \phi(y+u,\bary+v)=\sum_{k,l\geq 0} a_{k,l} u^k v^l .
\]
Setting $u=r\, e^{i\theta}$, $v=\baru$ we obtain for $r$ small enough
\[ 
   0 = \phi(y+u,y+\baru)= \sum_{m\geq 0}
        r^m \sum_{k=0}^m a_{k,m-k} e^{i (2k-m)\theta} ,
\]
which vanishes iff $a_{k,l}=0$ for all $k,l\geq 0$. \Halmos\\

 Consider the mirror extension, 
$\hatD=\DD\times \barDD\simeq \DD^2$,  of the unit disk, $\DD$.
We write $d_{\tinyD} = \frac{4 dz \bardz}{(1-z \barz)^2}$ for the Poincar\'e
metric on $\DD$. [By abuse of notation we write
$dz \bardz$ for the symmetric two tensor, 
$\half (dz\otimes \bardz + \bardz\otimes dz)$]. Also note, that
when $c$ is a complex number,
  $dz (c\frac{\partial}{\partial z}) = c$, but
  $dz (c\frac{\partial}{\partial \barz}) = \bar{c}$ (and not zero!).
Below, we will use the following metric on $\hatD$~:
\[
 d^{(1)}_\htinyD= 
      \frac{ |dz_1|}{1-z_1\barz_1}+
      \frac{|dz_2| }{1-z_2\barz_2} \equiv ds_1 + ds_2,
\ \ \ (z_1,z_2)\in \hatD.
\]
 This metric is more convenient here than the 
Riemannian metric,
$d^{(2)}_\htinyD= \sqrt{ds_1^2+ds_2^2}$.

\begin{Definition}
We denote by $\Aut(\DD)$ the group
 of holomorphic automorphisms of the disk consisting of
all M\"obius transformations which may be written
$R(z)=\frac{az+b}{\barb z+\bara }$, $|a|>|b|$. 
 To each $R\in \Aut(\DD)$ write
 $\barR(w)\equiv \overline{R(\barw)}$, $w\in \DD$ for the conjugated map.
The pair $\hatR=(R,\barR)$ acts isometrically on the extension,
 $\hatR^* d_\htinyD=d_\htinyD$, and
preserves the mirror diagonal. We denote by
 $\Aut(\hatD;\diag \DD)$ the collection of such pairs
and call it the group of mirror automorphisms of $\hatD$.
It is a subgroup 
of $\Aut(D^2)$ which itself has a fairly simple explicit description,
 see e.g.\ \cite[Proposition 11.1.3]{Kran00}.
\end{Definition}

\begin{Proposition}
\label{Proposition two form}
The holomorphic two-form $\hatg$ given by 
   \begin{equation}
      \hatg = \frac{4 dz_1 dz_2}{(1-z_1z_2)^2}.
     \end{equation}
 is the unique symmetric holomorphic two form on
 $\hatD=\DD\times \barDD =\DD^2$ which extends the Poincar\'e metric
 on the diagonal, i.e.\ such that
 \[
   d_{\mbox{}_{\Bbb D}} = j^* \hatg .
 \]
A mirror automorphism preserves
the holomorphic two-form, i.e.\  for all $\hatR \in \Aut(\hatD; \diag \DD)$,
  \begin{equation}
     \hatR^* g_\htinyD = g_\htinyD .
     \label{two form invariance}
   \end{equation}
\end{Proposition}
Proof: A calculation shows that indeed we obtain an extension. 
By the previous Lemma,
the factor $1/(1-z_1z_2)$ is uniquely determined by its value on the
diagonal. The assertion 
     (\ref{two form invariance})
is equivalent to the identity,
 \[ R'(z_1)\barR'(z_2) \frac{(1-z_1z_2)^2}{(1-R(z_1)\barR(z_2))^2} \equiv 1,
       \ \ \ \forall (z_1,z_2)\in \hatD ,\]
which is seen either by direct calculation or by the fact that it is
indeed correct on the mirror diagonal
 (where it expresses the fact that $R\in \Aut(\DD)$) and then by unicity of
mirror holomorphic functions,
Lemma \ref{unique KD}(3)
\Halmos\\

Let $\psi:\DD \rr \DD$ be a holomorphic map without critical points.
The pull-back of the Poincar\'e metric by $\psi$
 is proportional to the Poincar\'e
metric itself, where the factor of proportionality 
precisely defines the (square) of
the conformal derivative,
   \[ \psi^* d_\tinyD = (D\psi)^2 d_\tinyD, \ D\psi>0 .\]
It is independent of choice of conformal coordinates on $\DD$, i.e.\
under conjugations by $R\in \Aut(\DD)$ in either the source or in the
image. 
We write $\barpsi(z)\equiv \overline{\psi(\barz)}$, $z \in \DD$ for
the associated conjugated map.
The mirror extended map, $\hatpsi=(\psi,\barpsi)$,
is the unique map of $\hatD$ for which $\hatpsi \circ j = j \circ \psi$. It
preserves the diagonal 
but is, in general, not conformal on $\hatD$ 
(with respect to neither $d^{(1)}_\htinyD$ nor $d^{(2)}_\htinyD$).
It is, however, `conformal' with respect to our holomorphic two-form,
$g_\htinyD$. More generally, if
$\psi_1,\psi_2:\DD \rr \DD$ are two holomorphic maps, then their direct
product $\Psi=(\psi_1,\psi_2)$ verifies,
 \[     \Psi^* g_\htinyD = (D\Psi)^2 g_\htinyD ,\]
with a `conformal' derivative given by the formula,
\[
 (D\Psi)^2 = \psi_1'(z_1)\psi_2'(z_2)
                 \frac{(1-z_1 z_2)^2}{(1-\psi_1(z_1)\psi_2(z_2))^2} .
\]

Let $\hatpsi$ be the above mirror extension of $\psi$.
Then $(D\hatpsi)^2$ is
real and strictly positive on the mirror diagonal. We may then define
$D\hatpsi$ as the unique positive square root on the mirror 
diagonal and extend holomorphically to all of $\hatD$. On the mirror diagonal
it coincides with the usual definition of the conformal derivative
of $\psi$ on $\DD$, i.e.\
  \[
     j^* D\hatpsi^2 = D\psi^2 .
  \]
 Also, when $\Psi$ is a continuous
deformation of a mirror extended map,
$\Psi=\hatpsi$, then we may  define
$D\Psi$ by following the square-root along
the deformation (again provided that 
there are no critical points).

\begin{Lemma}
Let $\Psi=(\psi_1,\psi_2)$ be a direct product map on $\hatD$.
Then  for $i=1,2$,
  \[
       (1-|z_i|^2) \left|  
           \frac{\partial}{\partial z_i}
              \log   (D\Psi^2) (z_1,z_2) \right|
   \]
is conformally invariant with respect to 
mirror automorphisms, $\hatR \in \Aut(\hatD; \diag\DD)$,
both in the source and in the image.
\end{Lemma}
Proof: To see this, we consider maps $\hatR_1,\hatR_2\in\AutDD$
and the conjugated direct product,
$\Phi=\hatR_2\circ \Psi\circ \hatR_1$. 
Since $D\hatR_i^2\equiv 1$, $i=1,2$, we have that
$D\Phi^2 = D\Psi^2\circ \hatR_1$. Let 
$(z_1,z_2)=\hatR_1(u_1,u_2)=(R_1(u_1),\barR_1(u_2)$.
Taking the derivative with respect to $u_1$ and using
$|\partial R_1/\partial u_1| = (1-|z_1|^2)/(1-|u_1|^2)$ we obtain
\[
       (1-|z_i|^2) \left|  
           \frac{\partial}{\partial z_i}
              \log   (D\Psi^2) (z_1,z_2) \right| =
       (1-|u_i|^2) \left|  
           \frac{\partial}{\partial u_i}
              \log   (D\Phi^2) (u_1,u_2) \right|
\]
and thus the desired conformal invariance.\Halmos\\

Let $\psi:\DD\rr \DD$ be a conformal map without critical points.
\begin{Definition}
\label{def injectivity radius}
We define the {\em injectivity radius}
 $r=r[\psi](z) \in ]0,+\infty]$  
 of $\psi$ at $z\in \DD$ as
the largest value such that $\psi$ is injective on
a disc of hyperbolic radius $r$, centered at $z$. 
(In analogy with a similar
notion for Riemann surfaces, see e.g.\ \cite[section 2.9]{McMul94}).
 We call $\rho=\rho[\psi](z)=\tanh \frac{r}{2}$ the
{\em Euclidean radius of injectivity}.
If $R\in \Aut(\DD)$ maps  zero to $z$,
then  $\psi\circ R$ 
 is precisely injective on the Euclidean disc $B_\CCs(0,\rho)$.
\end{Definition}

\begin{Proposition} (Mirror K\"oebe distortion)
  \label{Prop mirror koebe}
   Let $\Psi=(\psi_1 ,\psi_2)$ be a direct product map on $\hatD$
   where both maps $\psi_1$ and $\psi_2$ are conformal maps from
   $\DD$ into itself and without critical points. At a given
   point $(z_1,z_2)\in\hatD$ we write $\rho_i=\rho[\psi_i](z_i)$,
    $i=1,2$ for the corresponding Euclidean radii of injectivity.
   We then have
        \begin{equation}
               \label{mirror Koebe distortion}
              | d \log \; D\Psi^2 | \leq 
                 (2+\frac{4}{\rho_1})ds_1 +
                 (2+\frac{4}{\rho_2})ds_2 . 
        \end{equation}
 \end{Proposition}
Proof:
We will use conformal invariance twice. 
  Let $(z_1,z_2)\in \hatD$. 
Fix mirror automorphisms
 for which
$\hatR_1(0,u_2)=(z_1,z_2)$ and $\hatR_2(\psi_1(z_1),\psi_2(z_2))=(w_1,0)$.
The conjugated product map,
$\Phi=\hatR_2 \circ
 \Psi \circ \hatR_1=(\phi_1,\phi_2)$, then maps $(0,u_2)$ to 
 $(\phi_1(u_1),\phi_2(0))=(w_1,0)$. The conformal derivative at
 $(u_1,u_2)$ is then given by
  \begin{equation}
  D\Phi^2 (u_1,u_2)= \phi_1'(u_1) \phi_2'(u_2) (1-u_1u_2)^2.
  \label{eq mirror derivative}
 \end{equation}
 
Therefore, 
 \begin{equation}
        (1-|u_1|^2)
      \frac{\partial}{\partial u_1} \; \log\; D\Phi^2(u_1,0)_{|u_1=0} = 
      \frac{\phi_1''(0)}{\phi_1'(0)} - 2u_2.
\label{intermediate Koebe}
\end{equation}

Here, $\phi_1$
is univalent on the 
Euclidean disk of radius 
$\rho_1=\rho[\psi_1](z_1)$ centered at zero. 
By the standard 
K\"oebe estimate, $|\phi_1''(0)/\phi_1'(0)|\leq 4/\rho_1$.
Also $|u_2|\leq 1$
 (in fact, $|u_2|=\tanh(d(z_1,\barz_2)/2)$ for a slightly better estimate).
The right hand side of 
(\ref{intermediate Koebe}) therefore does not exceed 
$(2+4/\rho_1)$.
Combining with the previous Lemma
on the  conformal invariance of the distortion, we obtain
\[  
          \left|  
           \frac{\partial}{\partial z_i}
              \log   (D\Psi^2) (z_1,z_2) \right| \leq
        (2+\frac{4}{\rho_1})\; \frac{1}{ 1-|z_i|^2} .
\]
 Noting
that $ds_1=|dz_1|/(1-|z_1|^2)$ and 
including the same estimate for the second variable  we obtain the
desired bound. \Halmos\\
      
Let us now return to our hyperbolic space, $U\subset \widehat{\CC}$
and the compact subset $K\subset U$.
We define the constants $\Delta$ and $\Delta'$ as in
(\ref{def Delta}).
 Let $K_\Delta=N_\Delta(K)\subset U$ be the $\Delta$
neighborhood of the compact set $K$. 
Below we will 
make use of constants,
\begin{equation}
   \ell_2(K,U,\Delta)>0\ \ \mbox{and}\ \
       \alpha_2(K,U,\Delta) >0,
    \label{def ell2}
\end{equation}
defined as follows: Consider $x\in K$, $u\in K_\Delta$ and let
$\gamma_{x,u}$ be a shortest geodesic between the two points.
 We let $\ell_{x,u} \in ]0,+\infty]$ be the
minimal length of non-contractible closed geodesics intersecting 
$\gamma_{x,u}$. Finally we let $\ell_2=\ell_2(K,U,\Delta)$ be the infimum
of all such lengths $\ell_{x,u}$. Because of $K$ and (the closure of)
$K_\Delta$ being compact sets, this infimum is necessarily non-zero.
We set $\alpha_2(K,U,\Delta) \equiv  \tanh(\ell_2/4) \in ]0,1]$ which is 
then a lower bound for the Euclidean radius of injectivity of 
a Riemann mapping  centered at a point along any of the above
mentioned geodesics.

\begin{Remarks}
In the previous section our definition of $\ell$ was somewhat simpler
because a local distortion estimate sufficed. Below we need a
global estimate and for this we need  to control the distortion along
paths (geodesics) connecting points in $K_\Delta$.
\end{Remarks}

\begin{Proposition}[Global Mirror Distortion]
\label{Global Koebe}
Let $f_1,f_2\in \EE(K,U)$ and consider in the cover,
 as in the previous section,
`inverses' 
$\psi_1:\DD\rr \DD$ and $\psi_2:\DD\rr \DD$ of $f_1$ and
$\bar{f_2}$, respectively. We write 
$\Psi=(\psi_1,\psi_2): \hatD \rr \hatD$ for the product map.
For $\xi\in\KD$ and $\zeta\in\diagK$.
we have the bound
(with the constant 
$\alpha_2$ from 
    equation (\ref{def ell2})),
\[
      \left| \log \frac{D\Psi^2(\xi)}{D\Psi^2(\zeta)} \right|
       \leq  (2+\frac{4}{\alpha_2}) \, d_{\hatU}(\xi,\zeta).
\]
\end{Proposition}
Proof:
Let $\gamma=(\gamma_1,\gamma_2)$  be a shortest geodesic between
 $\xi$ and $\eta$ for the metric, $d^{(1)}$, i.e.\  
the line element $ds=ds_1+ds_2$.
Then $\gamma_1$ and $\gamma_2$ are then shortest
geodesics between the two coordinate projections of $\xi$ and $\eta$,
Along $\gamma$ we have by Proposition 
  \ref{Prop mirror koebe} the
infinitesimal inequality 
 $|d\, \log D\Psi^2 | \leq (2 + \frac{4}{\alpha_2})\, ds$ and
the result follows by integration from $\log 1 = 0$.
\Halmos\\

The previous distortion Lemma is, in reality, only for the extended disk.
We need to establish distortion estimates for the mirror extension of $U$.
Unfortunately, it is not possible to do so globally
(unless $U$ is simply connected). We consider instead a restriction to
a neighborhood of the mirror diagonal.

\begin{figure}
\begin{center}
\epsfig{figure=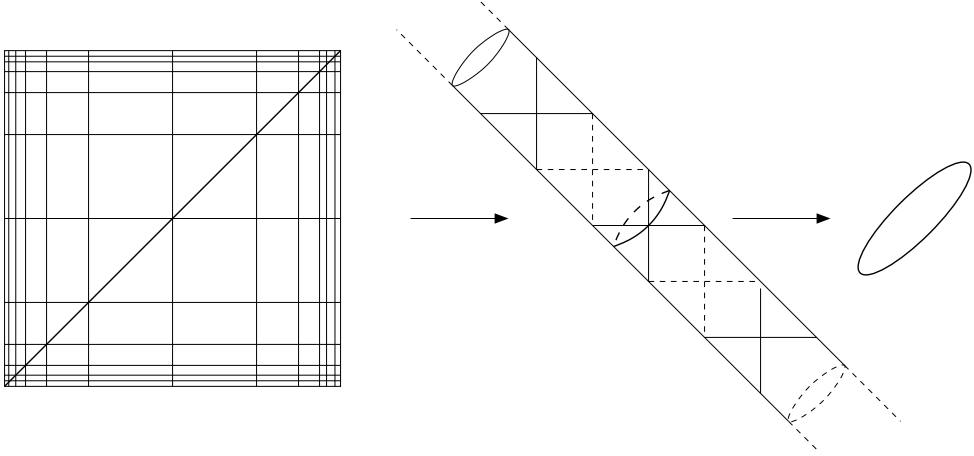,height=6cm}
\label{DmodK}
\end{center}
\caption{The quotients\~:
 $\hatD \rightarrow \hatD/\diag \Gamma \rightarrow \DD/\Gamma$.
 The illustration is in the case of an annulus where $\Gamma$ is generated
 by one element only. Since
 4 dimensions is difficult to illustrate we have only drawn real sections.}
\end{figure}

Let $\Gamma\subset \Aut(\DD)$ be a surface group of $U$
consisting of all automorphisms of $\DD$ that leaves
invariant a given Riemann mapping $\phi:\DD\rr U$.
The mirror extended surface group,
$\diag \Gamma=\{ (R,\barR): R\in \Gamma\}$ acts
`diagonally' upon $\hatD$ (it is a subgroup of $\Aut(\hatD;\diag K)$
and is normal iff $\Gamma$ is Abelian).
The quotient
(see Figure \ref{DmodK}),
   \[
      \hatD / \diag \Gamma,
   \]
is a complex  $2$-dimensional manifold.\footnote{
       $\hatD/\diag \Gamma$ could be viewed as a (non-trivial)
 fiber-bundle over $U$ with
 fiber $\barDD \simeq \DD$.}
Proposition
\ref{Proposition two form}
       shows that $g$ passes down to the
quotient as a holomorphic 2 form. The same is true
for the distortion estimate in
Proposition
\ref{Global Koebe}.

For $r>0$ denote by $\HK_r\equiv N_r(\diag K)$ the $r$-neighborhood
of $\diag K$ in $\hatU$. 
We lift $\HK_r$ to
     the set $\HK_{r,\Gamma}$ in $N_r(\diag\hatD /\diag \Gamma)$.
We claim that if 
$r<\ell/4$ then the natural projection,
  \[
          \HK_{r,\Gamma} \rr \HK_r,
  \]
is a conformal isomorphism (in particular, the lift consists of
une unique `copy' of $\HK_r$).
If this were not so then we could find
$z\neq z' \in \hatD/\diag \Gamma$ and
 $\eta,\eta'\in  \diag\hatD / \diag \Gamma$
for which $d(z',\eta)<r$, $d(z',\eta')<r$ and 
$\pi(z)=\pi(z')$, i.e.\ projects to the same point in $\HK_r$.
But then $d(\eta,\eta')<2r$ and there is a non-contractible
loop $\gamma=(\gamma_1,\gamma_2)$ containing e.g.\ $\eta$,
intersecting $\diag K$
and of length $\leq 4r < \ell$. Then at least one of $\gamma_1$ and 
$\gamma_2$ is non-contractible, of length $< \ell$ and intersects
$K$ and this is impossible.

Our two-form, $g$, on 
$\hatD/\diag \Gamma$ projects now to
a unique holomorphic two-form,
which we still denote $g$, on $\HK_{\ell/4}$.
This is the unique analytic continuation of the conformal metric
that we are searching for. It verifies,
    \[
       d_{\hatU | \diag U} = j^* g
   \]
Now, let 
    $\xi,v\in \HK_{\ell/4}$ and suppose that
$\Psi=(\psi_1,\psi_2): \OO(\xi) \rr \OO(v)$ 
is a locally defined product map between neighborhoods
(in $\HK_{\ell/4})$  of the two points.
We may then define the conformal
derivative of this map through the identification
\begin{equation}
     \Psi^* g_{|v} \equiv D\Psi^2(\xi) \; g_{|\xi}.
     \label{good conf deriv}
\end{equation}
When $\Psi$ preserves the diagonal, then 
$D \Psi^2_{| \diag K} >0$  and we may define its
positive  square root or
principal logarithm in the usual way.

% Not needed:
% If     $y \notin \diag \DD$    but   $\hatphi y \in \diag K$
% then its distance to    $\diag \DD$   is at least   $\ell$.
% Just take $\hatq=(q,\barq) \in \diag \DD$ and calculate
%    \ell \leq d(y_1,y_2) 
%         \leq d(y_1,q)+d(q,\bary_2)
%          =   d(y_1,q)+d(y_2,\barq) = d(y,\hatq)

\subsection{Mirror extended transfer operators and cone contractions}

Let $f\in\EE(K,U)$ and let $\hatf=(f,\barf)$ be the mirror extended map.
For $\eta\in \diag(K)$, we write for its mirror-preimages 
\[
   P_{\hatf}(\eta) \equiv 
      \hatf^{-1}(\eta)\cap \diag K \equiv \{u_i\}_{i\in J},
\]
 with $J$ an index set. 
We wish to define an analytic continuation of this ensemble
to points in $\KD$.
For $\xi\in \KD$, pick $\eta\in \diag K$ and a path
$\gamma$ in $\KD$ connecting $\eta \rr \xi$. 
For each $i\in J$,
$\gamma$ lifts by $\hatf$ to a path $\gamma_i$ connecting 
$u_i$ to some point $v_i\in \HK_{\Delta/\beta}$ (because of
contraction of the inverse map).  The collection
\[
P_{\hatf}(\xi)\equiv \{v_i\}_{i\in J} \subset \HK_{\Delta/\beta}
\]
yields the desired continuation.
 This set depends only on $\xi$ 
(and $f$, of course) but not on the choices of $\eta$ and
the  path $\gamma$.
Any other choice will just give rise to a permutation of $J$.
This is true if $\gamma$ is a shortest geodesic to the diagonal
(because its length is smaller than $\Delta<\ell/4$). But then it is
also true for any other path as long as the path stays within $\KD$.

Denote by $D\hatf^2(v)$, $v\in P_\hatf(\xi)$
 the holomorphic conformal derivative,
     (\ref{good conf deriv}),
 of $\hatf$.
 When $\eta$,
whence also $u\in P_\hatf(\eta)$,  belongs to the mirror
diagonal then $D\hatf^2(u)>0$ and we define its logarithm
by its principal value, $\log D\hatf (u) = \frac12 \log D\hatf^2(u)\in \RR$.
This extends to all $\xi\in \KD$, $v\in P_\hatf(\xi)$ by analytic
continuation, and arguing as above, is independent of the choices made.

Recall that
$A\equiv A(\KD)=C^0(\Cl\; \KD)\cap C^\omega(\KD)$ denotes the space 
of holomorphic functions having a continuous extension to the boundary.
We define
for $s\in \CC$, $\phi \in A(\KD)$ and $\xi\in\KD$ the transfer operator

\[
  L_{s,\hatf} \; \phi (\xi) 
     = 
 \sum_{v\in P_{\hatf}(\xi)}
      D\hatf(v)^{-s} \phi(v) \equiv
 \sum_{v\in P_{\hatf}(\xi)}
            e^{-s \log D\hatf(v)} \phi(v) .
\]

For the moment let us fix a real value of $s\geq 0$. Then 
$L_{s,\hatf}$ preserves $A_{\RRs}$,
the space of  self-adjoint elements.
We define for $\sigma>0$
a closed proper convex cone in $A_\RRs$,
 \[
   \CCC_\sigma = \{ \phi \in A_{\RRs}: |\phi(\xi)-\phi(\zeta)|
           \leq \phi(\zeta) (e^{\sigma d(\xi,\zeta)}-1),
	     \xi\in \KD, \zeta\in \diag{K} \}.
\]

We define $\beta(\phi_1,\phi_2)=
 \inf\{\lambda>0: \lambda\phi_1-\phi_2\in \CCC_\sigma \}$
and write $\dC=\frac12\;\log \; 
 (\beta(\phi_1,\phi_2) \beta(\phi_2,\phi_1))$
 for the corresponding projective Hilbert metric
(cf.\ \cite{Bir67,Liv95,Rugh02}).

Fix a point $\zeta_0=(x_0,\barx_0)\in\diag(K)$ and
denote by $\ell\in A_\RRs$ the (real-analytic) linear functional
\[
     \ell(\phi) = \phi(\zeta_0), \ \ \phi\in A.
\]
We use this to introduce  the {\em sliced cone},
\[
     \CCC_{\sigma,\ell=1} \equiv \{ \phi\in \CCC_\sigma : \ell(\phi)=1\}.
\]

\begin{figure}
\begin{center}
\epsfig{figure=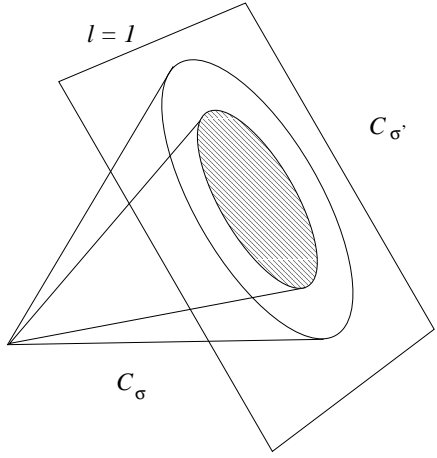,height=6cm}
\label{cones}
\end{center}
\caption{The cone contraction. 
The sliced cone $\CCC_{\sigma',\ell=1}$ has 
  an $R$-neighborhood which is contained in $\CCC_\sigma$.}
\end{figure}

\begin{Lemma}
[Cone contraction]
\label{Cone contraction}
Let $s\geq 0$ and choose $\sigma\equiv \sigma(s)>0$ large enough so that
\[
\sigma'\equiv \sigma'(s) =  \left(  1+ \frac{2}{\alpha_2}  \right) s + 
          \frac{1}{\beta} \, \sigma <\sigma.
\]
 Then there is $\eta<1$
such that for every $f\in \EE(K,U)$, the operator
 $L_{s,\hatf}$
 maps  $\CCC_\sigma$ into $\CCC_{\sigma'}$ and is an 
$\eta$-Lipschitz contraction for the Hilbert metric, $d_{\CCC_\sigma}$.
 Furthermore,
\begin{enumerate}
\itemsep -1mm
 \item[(a)] 
       There is $k>0$ such that  $\ell(\phi) \geq k \|\phi\|$
          for all $\phi\in \CCC_\sigma$ 
        (this corresponds to the dual cone having non-empty interior;
        we say that the cone is outer regular).
 \item[(b)] There is $R>0$ such that if 
          $\phi\in \CCC_{\sigma',\ell=1}$ then
	  $B(\phi,R)\subset \CCC_{\sigma}$ ($\CCC_\sigma$ has a
	    `uniformly large'
	  interior; we say that the cone is inner regular).
\end{enumerate}
\end{Lemma}
Proof:
Fix $\xi=(\xi_1,\xi_2)\in \KD$,
 $\eta=(\eta_1,\eta_2=\bar{\eta_1})\in \diag(K)$ 
and let $\hat{\gamma}=(\gamma_1,\gamma_2)$ be a shortest
geodesic for the metric $ds=ds_1+ds_2$ in $\hatD$.
Then $\gamma_i$ is a shortest geodesic between $\xi_i$ and $\eta_i$,
$i=1,2$. We write
 $d=d_{\hatU}(\xi,\eta)=\len(\gamma_1)+\len(\gamma_2)$ for the total
length. By considering pre-images by $F\equiv\hat{f}$ of
 $\hat{\gamma}$ we obtain
a pairing $(u,v)$ of the corresponding pre-images of $\xi$ and $\eta$
which verifies,
  \[
      d_{\hatU}(u,v)\leq \beta^{-1} d,
  \]
because of contraction of the local inverse maps. By definition 
of $\alpha_2=\alpha_2(K,U,\Delta)$ and the use of 
Mirror K\"oebe distortion, Proposition
\ref{Global Koebe}, we obtain,
\begin{eqnarray*}
   \lefteqn{\|
     L_{s,F} \phi(\xi)
     - L_{s,F} \phi(\eta) \|}\nonumber \\
         & \leq & \sum \left| 
	        (DF(u))^{-s} \phi(u) -
	        (DF(v))^{-s} \phi(v) \right| \\
  & \leq &  
         \sum |DF(u)|^{-s} |\phi(u) -\phi(v)| +
	        | DF(u)^{-s} - DF(v)|^{-s}|
		  \phi(v)  \\
  & \leq &  
        (e^{(1+\frac{2}{\alpha_2})s d} (e^{\beta^{-1} \sigma d} -1)
	+ (e^{(1+\frac{2}{\alpha_2})s d} -1) ) L\phi(\eta) \\
  & \leq &  
	(e^{{\sigma'} d} -1)  L\phi(\eta) \\
\end{eqnarray*}
where $\sigma'=(1+\frac{2}{\alpha_2})s+ \beta^{-1} \sigma$.

Thus, $L\equiv L_{s,\hatf}: \CCC_\sigma \rr \CCC_{\sigma'}$
and we get for the projective diameter
 (for this standard calculation we refer to e.g.\
 \cite{Liv95} or \cite[Appendix A]{Rugh02}),
 \[
   \diam_{\CCC_{\sigma}} \CCC_{\sigma'} \leq
         D = 2\log \frac{\sigma+\sigma'}{\sigma-\sigma'} 
         + \sigma' \diam \KD < \infty,
\]
 where we write $|K|$ for the diameter of $K$ in $U$.
 By Birkhoff's Theorem 
(see \cite{Bir67},\cite{Liv95} or \cite[Lemma A.4]{Rugh02})
  this implies a uniform contraction for the
hyperbolic metric on $\CCC_\sigma$. 
Writing $\eta=\tanh(D/4)<1$ we have for $\phi_1,\phi_2\in\CCC_\sigma$,
  \[ 
      \dC(L\phi_1,L\phi_2)\leq \eta \ \dC(\phi_1,\phi_2).
  \]

Property (a) is clear from the definition of the cone which implies~:
 \[
   |\phi(\xi)| \leq
            \ell(\phi) \;e^{\sigma \diam \KD}, \ \
          \phi\in\CCC_\sigma,\,
             \xi\in\KD.
\]

Set $\kappa = \frac{2}{\tanh(\Delta/2)}$ and
let $\phi\in A(\KD)$. We claim that for 
$\zeta\in\diagK$ and $\xi\in\KD$:
\begin{equation}
    |\phi(\xi)-\phi(\zeta)| \leq 
 |\phi| \, \kappa \, d(\zeta,\xi).
   \label{distortion phi}
\end{equation}

 It suffices to verify this inequality in
the universal cover. Consider coordinates
where $0\in \diag \DD \mapsto \eta$ and $u\in \hatD\mapsto \xi$.
Let $R=\tanh(\Delta/2) / (|u_1|+|u_2|)$. If $|t|\leq R$ then
$t|u_1|+t|u_2|\leq \tanh (\Delta/2)$ which implies
$d_\DD(tu_1,0)+d(tu_2,0)\leq \Delta$ and then also
$d^{(1)}(tu,0)\leq \Delta$. It follows that
$t\in B_\CCs(0,R) \mapsto \phi(tu)$ is analytic and bounded by $\|\phi\|$.
By the Schwarz Lemma we obtain 
\[
      |\phi(\xi)-\phi(\eta)| \leq 
          \frac{2 \|\phi\|}{R} \leq
          \frac{2 \|\phi\|}{\tanh (\Delta/2)} \, d(\xi,\eta) .
\]

Consider $h\in\CCC_{\sigma',\ell=1}$ and $\phi\in A_\RRs(\KD)$.
In order for $h+\phi$ to belong to $\CCC_\sigma$ we need that
\[
  \left| \frac{h_\xi+\phi_\xi}{h_\eta+\phi_\eta} - 1 \right|  \leq
        \exp({\sigma} d(\xi,\eta))-1
\]
is verified for all $\xi\in \KD$ and $\eta\in \diag(K)$. 
Setting $d=d(\xi,\eta)$ we see that
this is the case if 
\[
    |h(\xi)-h(\eta)|+|\phi_\xi-\phi_\eta|
\leq (h(\eta)-\|\phi\|) (e^{\sigma'd}-1).
\]
Using that $h\in \CCC_{\sigma'}$, $\ell(h)=1$ and the above 
distortion estimate
 (\ref{distortion phi}) for $\phi$  we see that it suffices that
for all $d>0$,
\[
     \|\phi\| \leq 
           \frac{e^{\sigma d}-e^{\sigma' d}}
                {\kappa d+e^{\sigma d} -1}
                  e^{-\sigma' \diam K}.
\]
Letting $d\rr 0$ the right hand side tends to
      $(\sigma-\sigma') \exp(-\sigma' \diam K) / (\kappa+\sigma) >0$
and in the $d\rr \infty$ limit it tends to $\exp(-\sigma' \diam K)>0$.
It follows that it has a minimum $R>0$ and we have shown property (b).
\Halmos\\

Consider now a sequence  $L_1,L_2,\ldots$ of operators
as in the above Lemma. 
We write $L^{(n)}=L_n \circ L_{n-1} \circ \cdots \circ L_1$ for the 
$n$'th iterated operator.
\begin{Lemma}
\label{cone bounds}
There are constants $c_1,c_2<\infty$ such that for
$h,h'\in \CC_{\sigma',\ell=1}$, $\phi\in A$ and $n\geq 1$:
\begin{enumerate}
  \item[(1)] \ \ \
	    $\displaystyle
            \left| \frac
               {L^{(n)}h}{\ell(L^{(n)}h)} -
               \frac{L^{(n)}h'}{\ell(L^{(n)}h')} \right|
             \leq c_1 \eta^n
            $,
  \item[(2)] \ \ \
	    $\displaystyle
            \left| \frac
               {L^{(n)}\phi}{\ell(L^{(n)}h)} -
               \frac{L^{(n)}h}{\ell(L^{(n)}h)}\;
               \frac{\ell(L^{(n)}\phi)}{\ell(L^{(n)}h)}
               \right|
             \leq c_2 \eta^n |\phi|
            $.
\end{enumerate}
\end{Lemma}
Proof: Outer regularity, i.e.\
Property (a) of the above Lemma, and a computation show that
for $\phi_1,\phi_2\in\CCC_{\sigma, \ell=1}$,
  \[
      \label{k bound}
      |\phi_1-\phi_2| \leq \frac{k+1}{k^2} (e^{d(\phi_1,\phi_2)/2}-1).
  \]
When $\phi_1,\phi_2\in L^{(n)} \CCC_\sigma$, $\ell(\phi_1)=\ell(\phi_2)=1$
and $n\geq 1$ we have that
 $d_{\CCC_{\sigma'}}(\phi_1,\phi_2)\leq \eta^n D$ and therefore,
 \[
    |\phi_1-\phi_2| \leq \frac{k+1}{k^2} (e^{\eta^n D} -1 )
  \]
which is smaller that $c_1 \eta^n$ for a suitable choice of $c_1$.
This yields the first bound.

For the second bound note that $B(h,R)\subset \CCC_\sigma$.
 For $\phi\in A_\RRs$ (small) and $h\in \CCC_{\sigma'}$,
$d_\CCC(h+\phi,h)\leq \frac{1}{R} |\phi| + o(|\phi|)$ and therefore 
\[
  \left| 
           \frac
              {L^{(n)} (h+\phi)}
              {\ell (L^{(n)} (h+\phi))}  -
           \frac
              {L^{(n)} (h)}
              {\ell (L^{(n)} (h))}
	      \right| \leq
         \frac{k+1}{k^2} \eta^n \frac{1}{R} |\phi| + o(|\phi|)
\]
By linearizing this bound 
(and loosing a factor of at most $\sqrt{2}$) we may extend
this bound to any complex $\phi\in A$ to obtain the second inequality
with $c_2 = {\sqrt{2}}\frac{k+1}{k^2} \frac{1}{R}$.
\Halmos\\

\subsection{Analytic conformal families and mirror extensions}
Let $\OO_\RRs\subset\RR^n$ be an open set containing the origin and let
 $\OO_\CCs\subset \CC^n$ be an open convex neighborhood, invariant
under complex conjugation.
Also, let
  $t\in \OO_\CCs \subset \CC^n \rr f_t \in \EE(K,U)$ be a continuous map.
\begin{Definition} 
  \mbox{}
  \label{Analytic family}
  \begin{enumerate}
  \item
  $(f_t)_{t\in \OO_\CCs}$
  is called an analytic family, if  the map
  $\{(t,z): t\in \OO_\CCs, z\in \D_{f_t}\} \mapsto f_t(z)\in\CC$
  is analytic.
  \item
  We say that the family $f_t$ verifies an $L$-Lipschitz condition
  (with $0< L < +\infty$) if for any  $z\in K_\Delta$,
    and any choice of local inverse $f_0^{-1}(z)$,
  the map $t\in\OO_\CCs \mapsto \log Df_t \circ f_t^{-1}(z)\in\CC$ 
   is  $L$-Lipschitz.
  \item We define the condition number of $f\in\EE(K,U)$ to be
        \[
           \Gamma(f) = 
               \|Df \|_{f^{-1}K} \,
               \|1/Df \|_{f^{-1}K}.
       \]
  \end{enumerate}
\end{Definition}

  It is no lack of generality to assume that the parameters
are one-dimensional (n=1). We may also  assume that
 $\OO_\CCs=\DD$, i.e.\ is the unit-disk and consider 
$\OO_\RRs= \DD\cap \RR=]-1,1[$ as a real section.
In the following let
 $t\in\DD\mapsto f_t\in \EE(K,U)$
be an analytic family,
 verifying an
$L$-Lipschitz condition.

\begin{Notation}
\label{notation deriv}
 Below it is convenient to introduce
\[
\frac{d}{Dt}=(1-t\bart) \frac{\partial}{\partial t}
\]
 for the
conformal derivative from $(\DD,d_\DDs)$ to $\CC$ (with the Euclidean
metric). For a holomorphic map,
$h:\DD\rr \DD$, from the disk to itself
 we write also
\[
\frac{D}{Dt}=\frac{1-t\,\bart}{1-h(t)\,\overline{h(t)}}
 \frac{\partial h}{\partial t}
\]
for the conformal derivative between the disks. (Note that we do not
take absolute values).
\end{Notation}

Let $\hatU=U\times \barU$ be the mirror extension of $U$.
Let $\Delta>0$ be chosen as in the previous section.
We denote by $d_{\hatU}$
the metric on the mirror extension of $U$, 
induced by the metric
 on the universal cover, $\hatD$.
We obtain a conjugated analytic family if we set
$\D_{\barf_t}\equiv\overline{(\D_{f_\bart})}\subset \barU$ and for
$x'\in \D_{\barf_t}$, $\barf_t(x')\equiv \overline{f_{\bart}(\barx')}$.
Then $\barf_t(x')$ is analytic in $t$ and $x'$ on
$\{(t,z): t\in \DD, z\in \D_{\barf_t}\}$.
 We also define for $t\in\DD$  the product map
$F_t : (x,x')\in \D_{f_t} \times \D_{f_\bart}
\mapsto  (f_t(x),\barf_t(x'))\in\UbarU$.
Again, this map is analytic in $x,x'$ and $t$ on its domain of definition.

Consider
 $x_0\in \D_{f_0}$, $y_0=f_0(x_0)\in K$ and 
 Riemann mappings $\phi_{x_0}$ and $\phi_{y_0}$ defined as above.
Let $\psi_0$ be the corresponding `inverse' map.
Since $\DD$ is simply connected there is a unique holomorphic extension
\[
       (t,z)\in \DD\times \DD \mapsto \psi_t(z)\in \DD
\]
which analytically continues $\psi_0$ and 
defines an inverse (branch) of $f_t$, $t\in\DD$.
For fixed $z\in\DD$ the map $t\in\DD\mapsto \psi_t(z)\in\DD$
is holomorphic, whence it is non-expanding
and its conformal derivative can not exceed one, i.e.\
 \[
      |\frac{D}{Dt} \psi_t (z)| \equiv  
        \frac{1-t\bart}{1-\psi_t\overline{\psi_t}} | \partial_t \psi_t|
 \leq 1
  \]

The map $\barpsi_t(w)\equiv \overline{\psi_\bart(\overline{w})}$
defines an inverse of $\barf_t$ (in the corresponding cover). 
Also, $\Psi_t (w_1,w_2)=(\psi_t(w_1),\barpsi_t(w_2))$ defines
an inverse of $F_t$ in the mirror-extended cover of $U$.
Note that for $t$ real the functions $\psi_t$ and $\barpsi_t$ are
complex  conjugated and $\Psi_t$ therefore preserves the mirror diagonal
but that this is no longer true when $t$ becomes complex.\\

\begin{figure}
\begin{center}
\epsfig{figure=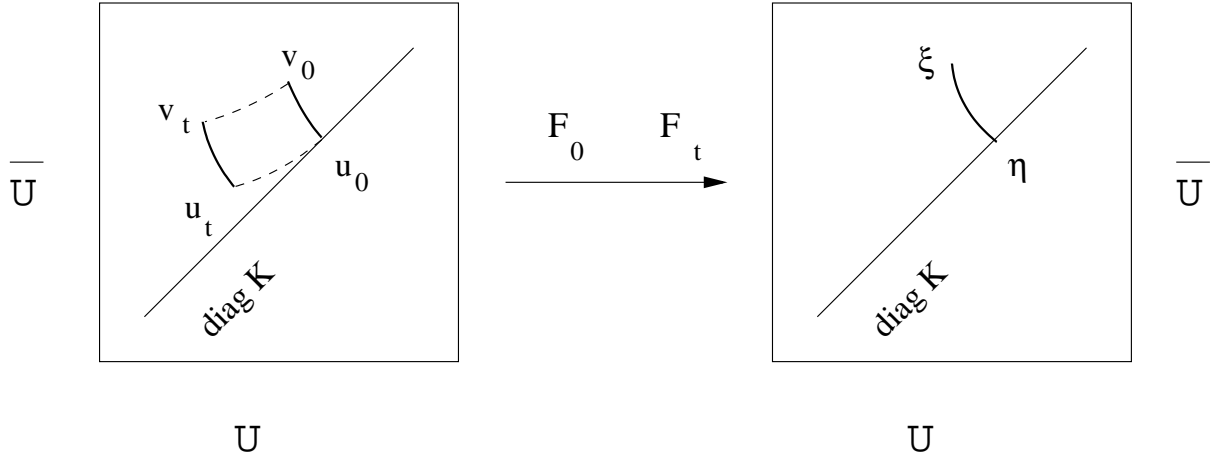,height=6cm}
\label{Finverse}
\end{center}
\caption{An inverse $v_t\in F_t^{-1}(\xi)$
 obtained by lift and analytic continuation.}
\end{figure}

For all $\xi\in\KD$, $t\in \DD$ we may define 
\begin{equation}
    P_{F_t}(\xi)=\{v^i_t\}_{i\in J}
\end{equation}
as follows: Let $\eta\in \diag K$ be at a distance $\leq \Delta$ to $\xi$.
Denote by $\gamma$ a shortest geodesic between the two points.
For every $u^i_0 \in \hatf_0^{-1} (\eta)\cap \diag K$, we lift $\gamma$
by $F_0$ to a path connecting $u^i_0$
to a point $v^i_0 \in F_0^{-1}(\xi)$ and we define then
$u^i_t\in F^{-1}_t(\eta)$ and also $v^i_t\in F_t^{-1}$ by analytic
continuation in $t\in \DD$ (by using a suitable inverse $\Psi_t$ in the
cover and projecting down). This defines the $v^i_t$ uniquely
up to permutations of $J$ (see Figure \ref{Finverse}).

If $\xi \in \hatU$ then
\[
   d(\Psi_t(\xi), \Psi_0(\xi)) =
 d(\psi_t(z_1), \psi_0(z_1)) +
 d(\barpsi_t(z_2), \barpsi_0(z_2))\leq
 2 d_\DDs(t,0) . 
\]
For $t$ real we know that
$\Psi_t: \HK_{\Delta} \rr \HK_{\Delta/\beta}$.
When making $t$ complex we want still to have a contraction of $\HK_\Delta$
and by the above it suffices to have
$\Delta/\beta+2 d_\DD(t,0) < \Delta$ or, equivalently,
  \begin{equation}
   |t| < \tanh  \left(\frac{\Delta}{4}(1-\frac{1}{\beta})\right).
     \label{first t condition}
  \end{equation}
When this condition is fulfilled we may analytically continue
the transfer operator in $t$.
First note that the conformal derivative,
 $DF_t(v)$  of $F_t$,
at a point $v\in P_{F_t}(\xi)$ is uniquely defined because
$v,\xi\in \KD\subset \widehat{K}_{\ell/4}$,
cf.\ equation (\ref{good conf deriv}). 
Recalling that $DF_0(u) >0$ when $u\in P_{F_0}(\eta)$, $\eta\in \diag K$
we may define $\log DF_0(u) \in \RR$ by its principal value and 
then $\log DF_t(v)\in\CC$ using
the lift of the path $\gamma$ and analytic
continuation in $t$.

 For $\phi\in A(\KD)$, $\xi\in \KD$, $s\in\CC$
and $t$ verifying (\ref{first t condition}), we set
  \[
   L_{s,F_t} \phi(\xi) =
       \sum_{v\in P_{F_t}(\xi)} DF_t^{-s}(v) \phi(v) 
       = \sum_{v\in P_{F_t}(\xi)} \exp( -s \, \log \, DF_t(v)) \phi(v) 
           .
 \]
This uniquely defines 
a bounded linear operator on $A(\KD)$.

The $L$-Lipschitz condition on $f_t$ is equivalent to the assumption that
$t\in (\DD,d_\DDs) \mapsto \log \,  D\psi_t$ is $L$-Lipschitz.
Since the maps are analytic we arrive at the equivalent condition,
\begin{equation}
 \left|  \frac{d}{Dt} \,\log \, D\psi_t \right| \leq L.
\label{L-log-Lip}
\end{equation}

Our hypotheses ensures that the local inverse of $F_t$,
i.e.\ the couple
$\Psi_t=(\psi_t,\barpsi_t)$, satisfies the conditions for the following
Lemma to apply~:

\begin{Lemma}
[Parameter distortion]
\label{Parameter distortion}
Let $\phi_{1t},\phi_{2t}$, $t\in \DD$ be holomorphic families of
conformal maps from the disk to itself, both having a conformal derivative
which is $L$-log-Lipschitz as in 
(\ref{L-log-Lip}).  Let 
$\Phi_t=(\phi_{1t},\phi_{2t}):\hatD\rr\hatD$ be their direct product.
Then
  \[
       \left|  \frac{d}{Dt} \log D \Phi^2 \right| \leq L+4
   \]
\end{Lemma}
Proof:
We have
(cf.\ the Notation \ref{notation deriv}),
 \[ D\phi^2 = |D_u\phi^2|=
    \partial_u \phi
    \overline{\partial_u \phi}
    \frac{(1-u\baru)^2}{(1-\phi \barphi)^2} .\]
Taking a derivative in $t$,
\[ \frac{d}{Dt} \log |D_u\phi^2| =
     ( \partial_u \phi )^{-1}
    \frac{d}{Dt}\partial_u \phi 
      + 
     2 \overline{\phi}\, {\frac{D}{Dt}} \phi.\]
In the identity,
\[  d \, \log D\phi^2 = 2 \,\Re \, \left( 
       \frac{d}{Dt} \log D\phi^2 \, \frac{dt}{1-t \bart}\right) ,
       \]
the left hand side is, by assumption,
 bounded in absolute value by $L |dt|/(1-t\bart)$.
But then $|\frac{d}{Dt} \log D\phi^2| \leq \frac{L}{2}$ and also,
\[
      \left|
(\partial_u \phi )^{-1}
\frac{d}{Dt} \partial_u \phi  
\right| \leq \frac{L}{2}+
            2 | \overline{\phi} {\frac{D}{Dt}}  \phi | \leq 
	    \frac{L}{2} + 2.
\]
Consider now the product map $\Phi=\Phi_t$.
Assume first that only $\phi_1$ depends on $t$. We then use the
same conjugation as in Proposition 
  \ref{Prop mirror koebe} 
to obtain the expression
  (\ref{eq mirror derivative}) for the conformal derivative $D\Phi^2$.
Taking now a $t$-derivative we get
  \[
    | \frac{d}{Dt} \log D\Phi^2| = 
  |
       (    {\partial_{u_1} \phi_{1t}} )^{-1}
      \frac{D}{Dt} \partial_{u_1} \phi_{1t}
| \leq \frac{L}{2}+2.
 \]
Adding the same contribution from the $t$-dependence of $\phi_{2t}$
we reach the desired conclusion.\Halmos\\

\begin{Lemma}
   Let $h\in \CCC_{\sigma'}$.
  Choose $x_0\in f_0^{-1} K$ and set $\lambda=\|Df\|_{f^{-1}K}$.
  Suppose that $d_\DD(t,0)\leq \Delta(1-\frac{1}{\beta})$.
  Let $\xi\in\KD$, $\eta\in \diag K$ and let
  $v_t\in P_{F_t}(\xi)$ and $u_t\in P_{F_0}(\eta)$ be pairs of pre-images
  constructed as above. Then
\[
  \left| 
  \frac
     {h(v_t) e^{-s \, \log \, (DF_t(v_t)/\lambda)}}
     {h(u_0) e^{-s_0 \, \log \, (DF_0(u_0)/\lambda)}}  - 1 \right|
      \leq e^q -1
\]
with 
\[
   q=|s-s_0| \log \Gamma (f_0) + 
   \left(2+|s|(2+\frac{L}{2})\right) d_\DD(t,0) +
   \left(\frac{1}{\beta} + |s| (1+ \frac{2}{\alpha_2})\right) d(\xi,\eta).
    \]
\end{Lemma}
Proof:
Since $h\in \CCC_{\sigma'}$ we know that
  \[ \left|\frac{h(v_t)}{h(u_0)}-1\right| \leq e^{d(v_t,u_0)}-1 .\]
and the distance in the exponent may be bounded as follows:
\[ d(v_t,u_0)\leq d(v_t,v_0)+d(v_0,u_0)
    \leq 2 d_\DD(t,0)+\frac{1}{\beta}d(\xi,\eta).\]
Lemma \ref{Global Koebe} and 
\ref{Parameter distortion} apply here so we 
also have the following inequalities,
    \[ \left|  \log  \frac{DF_t(v_t)}{DF_0(v_0)} \right| \leq
           \left(2+\frac{L}{2} \right) d_\DD(t,0) \]
    and
    \[ \left|  \log  \frac{DF_0(v_0)}{DF_0(u_0)} \right| \leq
           \left(1+\frac{2}{\alpha_2} \right) d(\xi,\eta).
    \]
By definition \ref{Analytic family} (3)
of the condition number of $f$,
 \[    \left|\log \frac{DF_0(u_0)}{\lambda}  \right| \leq \log \Gamma(f). \]
The inequality,
  \[   \left|  \prod e^{a_i} -1 \right| \leq e^{\sum |a_i|} -1 ,
   \]
is valid for any complex numbers, $a_1,\ldots,a_n$.
Now, insert the 4 estimates above  to obtain the
claimed inequality.\Halmos \\

The following non-linear map,
 \begin{equation}
    \pi_{s,F_t}(\phi) = 
        \frac{ L_{s,F_t} \phi }
        { \ell(L_{s,F_t} \phi) }
  \end{equation} 
is well-defined when the denominator does not vanish.

\begin{Lemma}[lemma neighborhood]
\label{lemma neighborhood}
Let $f_t\in \EE(K,U)$, $t\in \DD$ verify an $L$-Lipschitz
condition. For  $s_0\geq 0$ we let 
$W^{s_0}$
denote the open neighborhood of
 $(s_0,0)
\subset \CC\times \CC$
 consisting of all $(s,t)$ that verify
\[
   \left|s-s_0\right| \log \Gamma (f_0) + 
   \left(2+|s|(2+\frac{L}{2})\right) d_\DD(t,0) < \log\frac43,
    \]
and we let 
\[
   \rho  = \frac14  e^{-\sigma' \diam \KD} .
\]
Then for $h\in \CCC_{\sigma', \ell=1}$ and all $(s,t)\in W^{s_0}$,
$\phi\in A(\KD)$, $|\phi|<\rho$ we have
 \[
 1 \leq   \| \pi_{s,F_t} (h+\phi) \| \leq 5 e^{2 \sigma' \diam \KD}
\]
and also,
\begin{equation}
   \label{cone bound}
   \left| 
       \frac{\ell(\lambda^s L_{s,F_t} (h+\phi))}
            {\ell(\lambda^{s_0} L_{s_0,F_0} (h))}-1
   \right| \leq \frac23.
\end{equation}
\end{Lemma}
Proof: We first use our previous Lemma for $\xi=\eta$. We let $q$
and $u_t$, $v_t$ be as in that Lemma.
Our assumptions imply $e^q-1 < \frac13$ and therefore,
\[ 
\left|h(u_t)e^{-s \log \frac{DF_t(u_t)}{\lambda}} -
h(u_0)e^{-s_0 \log \frac{DF_0(u_0)}{\lambda}}
\right| \leq
\frac13 h(u_0)e^{-s_0 \log \frac{DF_0(u_0)}{\lambda}} .
\]
Summing this inequality over all
pairs of pre-images and then dividing by the right hand side,
we obtain
\begin{equation}
 \left|\frac{\ell(\lambda^s L_{s,F_t} h)}
   {\ell(\lambda^{s_0} L_{s_0,F_0} h)} - 1 \right| \leq
 \frac13.
   \label{first cone bound}
\end{equation}
In particular,
       $| \lambda^s \ell (L_{s,F_t}h) | \geq \frac23 
          \lambda^{s_0} \ell( L_{s_0,F_0} h)$,
for all $s$ and $t$ verifying the first condition.
Using the Lemma once more, and for general $\xi$ and $\eta$,
we see that
\[ 
|e^{-s \log \frac{DF_t(v_t)}{\lambda}}| \leq
\frac43 e^{-s_0 \log \frac{DF_0(u_0)}{\lambda}} \leq
\frac43 e^{\kappa} 
       e^{-s_0 \log \frac{DF_0(u_0)}{\lambda}}  h(u_0),
\]
where we have set $\kappa=\sigma' \diam \KD$. From this we obtain
\[ \| \lambda^s L_{s,F_t}\| \leq
       \frac43 e^{\kappa} \lambda^{s_0} L_{s_0,F_0}h(\eta) .
\]
When $\phi$ is of norm smaller than $\rho$ we have 
because $h\in \CCC_{\sigma',\ell=1}$,
\begin{equation}
    \label{second cone bound}
    \frac{|\ell(\lambda^s L_{s,F_t} \phi)|}
       {\ell(\lambda^{s_0}L_{s_0,F_0} h) }
     \leq \frac43 e^\kappa \rho=
         \frac13 .
\end{equation}
Finally, as $\|h\| \leq e^\kappa$, we obtain the upper bound~:
\[  \| \pi_{s,F_t}(h+\phi) \| =
         \frac{\| \lambda^s L_{s,F_t}(h+\phi)\|}
         {| \lambda^s \ell(L_{s,F_t}(h+\phi)) |}\leq 
     \frac{4/3 e^\kappa \|h+\phi\| }
          {2/3-1/3}
\leq 5 e^{2 \kappa}.
\]
 The lower bound is clear. 
The bound (\ref{cone bound}) follows from
(\ref{first cone bound})
and (\ref{second cone bound}).
\Halmos\\

\subsection{Analytic measurable sections}
Let us now return to the probability space $(\Omega,\mu)$ and
a $\mu$-ergodic transformation $\tau:\Omega\rr\Omega$.

We view the space $\Omega\times A$ as a (trivial) fiber
bundle over $\Omega$ with each fiber being $A=A(\KD)$. We denote by 
$\calA$ the set of measurable sections of this fiber bundle and
write $\|\Phi\|$ for the $\mu$-essential sup 
of an element $\Phi\in\calA$. $A$ is separable, so measurability and
Bochner-measurability is here the same.
Then $\calA$ is again a unital Banach algebra when we define the
analytic operations to be performed fiber-wise. We note that
measurability is preserved under such operations and also by
taking uniform limits. We write  $\calA_\RRs$ for the subspace of
real-analytic sections.
Let $\CCC_\sigma(\Omega)$ denote the space of measurable cone-sections
of $\Omega\times \CCC_\sigma$. 
We write $\CCC=\CCC_{\sigma,\ell=1}(\Omega)$ for the `sliced' 
measurable cone-sections. The latter forms a bounded subset of $\calA_\RRs$.

\begin{Assumption} Let $\OO\subset \RR^n$ be an open set and
  let $\OO_\CCs\subset \CC^n$ be a complex neighborhood of $\OO$.
  In the following we will assume that
  $t\in \OO \rr \FF^t = (\bff_{t,\omega})_{\omega\in\Omega}
  \in \EE_\Omega(K,U)$ is a map for which the following hold:
  \label{FF L Lipschitz}
  \begin{enumerate}
  \item
  For each $\omega\in\Omega$ the map
  $t\in\OO_\CCs\rr \bff_{t,\omega}\in\EE(K,U)$ is analytic in the sense
  of Definition \ref{Analytic family}.
  (Note that we are implicitly assuming
  that for each fixed $t\in\OO_\CCs$, the mapping 
  $\omega\in\Omega\mapsto \bff_{t,\omega}\in \EE(K,U)$
  is measurable as in Definition \ref{measurable family}).
  \item
  For each $t\in\DD$ and $\omega\in\Omega$, the map
   $\bff_{t,\omega}$, $t\in \DD$ verifies the $L$-Lipschitz condition
  in Definition \ref{Analytic family} for the same number $0<L< \infty$.
  \item The condition numbers $\Gamma(\bff_{0,\omega})$, $\omega\in\Omega$
    are uniformly bounded by some $\Gamma<+\infty$.
  \end{enumerate}
  \label{Assump}
\end{Assumption}

In the following we consider an analytic family,
$t\in\DD\mapsto \FF^t=(\bff_{t,\omega})_{\omega\in\Omega}\in \EE_\Omega(K,U)$
verifying Assumption \ref{Assump} above.\\

Let
$F_{t,\omega}=(\bff_{t,\omega},\overline{\bff}_{t,\omega})$ denote the
holomorphic extension  of $\bff_{t,\omega}$
and let $s_0=\dimH(J(\bff_{0,\cdot}))\in [0,2]$ be the
(a.s.) Hausdorff dimension
of the random Julia set at $t=0$.
We choose $\sigma=\sigma(s_0)$ and $\sigma'=\sigma'(s_0)$ so as to verify 
the Cone contraction conditions in Lemma \ref{Cone contraction}.
Let $W^{s_0}\subset \CC^2$ and $\rho>0$ be chosen as in
Lemma
\ref{lemma neighborhood}
and let $\bfh\in\CCC_{\sigma'}(\Omega)$. For 
$(s,t)\in W^{s_0}$ the following
 `sliced' cone-map,
\[
    \pi_{s,t}(\Phi)_\omega \equiv \pi_{s,F_{t,\omega}}(\Phi_{\tau \omega}) =
        \frac{L_{s,F_{t,\omega}} \Phi_{\tau \omega}}
        {\ell(L_{s,F_{t,\omega}} \Phi_{\tau \omega})},
\]
is a well-defined map $\pi_{s,t} : B(\bfh,\rho) \rr \calA$. By Proposition
\ref{lemma neighborhood} the image is bounded in norm by
$5 \exp( 2 \sigma' \, \diam \KD)$.
 It takes the value
of $\Phi$ at the shifted fiber $\tau \omega$, acts with the transfer
operator, normalises according to $\ell$ and assigns it to the fiber
at $\omega$. 
Measurability of the image is a consequence of the map 
$(s,F_t)\mapsto L_{s,F_t}$ being continous and
$\ell$ being strictly positive on the image. 
The reader may note that the (non-normalised)  family 
$(L_{s,F_{t,\omega}})_{\omega\in\Omega}$ need not be
uniformly norm-bounded, whence need not
define a bounded linear operator when acting
upon sections in $\calA$. This is 
the case e.g.\  in our example
in the introduction.

We denote by 
$\pi^{(n)}_{s_0,0}:
 \CCC_{\sigma'}(\Omega)\rr
 \CCC_{\sigma'}(\Omega)
 $ the $n$'th iterated map of
$\pi_{s_0,0}$ restricted to the cone-section.

\begin{Lemma}
\label{pi contraction} There are constants, $c_1,c_2<+\infty$ such that

\begin{enumerate}
\item  For $\bfh,\bfh'\in \CCC_{\sigma'}(\Omega)$ we have:
\[
    | \pi^{(n)}_{s_0,0}(\bfh) - \pi^{(n)}_{s_0,0}(\bfh')| \leq c_1\eta^n.
    \]
\item Taking the derivative of the $n$'th iterated map at the point
 $h\in\CCC_{\sigma'}(\Omega)$ we have
 \[
    \| D_\bfh \pi^{(n)}_{s_0,0}(\bfh)\| \leq c_2 \eta^n.
    \]
\item  The map, 
  \[
      (s,t) \in W^{s_0},\ \Phi\in B(\bfh,\rho)
        \mapsto \pi_{s,t}(\Phi) \in \calA
  \]
  is real-analytic.
\end{enumerate}
\end{Lemma}
Proof: (1) and (2) are reformulations of the bounds already given
in Lemma \ref{cone bounds} (with the constants from that lemma).  
A calculation shows that
\[
   (L_{s,F_{t,\omega}} \phi)^* = L_{\bar{s},\bar{F}_{t,\omega}}
       (\phi)^*, \ \ \phi\in A,
\]
which implies that for $s$ and $t$ real, the operator $L_{s,F_{t,\omega}}$
maps $A_\RRs$ into $A_\RRs$, i.e., is real-analytic.
 Each $\pi_{s,t}(\phi)_\omega$ is analytic in $s$, $t$ and $\phi$
 (for fixed $\omega$). Uniform boundedness was already shown 
above and
a Cauchy formula 
(choosing $r>0$ small enough),
    \[
       t \mapsto \left( \oint_{|t-t'|=r}
             \frac{\pi_{s,t}(\phi)_\omega}{t-t'}
              \frac{dt'}{2\pi i} \right)_{\omega\in \Omega}
    \]
enables us to recover a power series in the $t$-variable
(similarly for $s$ and $\phi$) within $\calA$. The map is real-analytic
in the sense that  it maps $(s,t)\in W^{s_0}\cap \RR^2$,
 $\Phi\in B(h,\rho)\cap \calA_\RRs$ into $\calA_\RRs$.\Halmos\\

First, we consider the real case 
$(s,t)\in W^{s_0}_\RRs \equiv W^{s_0}\cap \RR^2$.
 Let $h^0\equiv \bfone\in 
\subset \calA_\RRs$ be the
unit section of our bundle and define recursively the iterates 
$h^{k+1}=\pi_{s,0}(h^k)\in
 {\CCC_{\sigma',\ell=1}(\Omega)}$,
$k\geq 0$.
Lemma
\ref{cone bounds}
shows that $|h^{k+n}-h^k|\leq c_1 \eta(s)^k$ which tends exponentially
fast to zero.
The sequence thus converges uniformly in $\calA_\RRs$ towards a fixed point
  \[
       h^*=\pi(h^*) \in \CCC_{\sigma', \ell=1}(\Omega) .
  \]
We are interested in the normalisation factor, 
\[
 p_{s,t,\omega}=\ell(L_{s,F_{t,\omega}} h^*_{\tau\omega} )
 \]
at the fixed point. This function
	is real and strictly positive.

\begin{Lemma}
We have for  $s$ and $t$ real the following formula for the pressure~:
\[
   {P(s,\Lambda(\FF^t))} = \int \log p_{s,t,\omega} d\mu (\omega).
\]
\end{Lemma}

Proof: The embedding $j:K \rr \diag K \subset \hatU$ induces a pull-back
$j^* : \CCC_{\sigma'}\rr C(K)$. On $C(K)$ (before the mirror embedding)
we act with
the operator, $L_{s,f_{t,\omega}}$, as in 
(\ref{transfer op time}),
 and on the cone with the mirror
extended operator, $\LL_{s,\hat{\bff}_{t,\omega}}$.
Then  $L j^* h = j^* \LL h$ for $h\in \CCC_{\sigma'}$ and the
cone properties show that 
$\ell(h)\leq \|j^* h\| \leq \ell(h) e^{\sigma' \diam \KD}$.
  It then follows that
$\mu$-almost surely
\[
    P(s,\Lambda(\FF^t))=
    \lim_n \frac1n \log \|L^{(n)}_{s,{t,\omega}}\| =
    \lim_n \frac1n \log \ell(\LL^{(n)}_{s,{t,\omega}} h^*) =
     \lim_n \sum_k \frac1n  \log p_{s,t,\tau^k\omega}.
\]
The latter function is comparable to $\log \|L_{s,\bff_{t,\omega}}\|$,
whence integrable, so
by Birkhoff's Theorem it converges $\mu$-almost surely
towards the integral of $\log p$ as we
wanted to show.
\Halmos\\
 
\begin{Remarks}
The pressure does not depend on the choice of $\,\ell$ (of course, it
should not). If one makes another choice $\tilde{\ell}$ for the normalisation
this simply introduces 
a co-cycle that vanishes upon integration.
\end{Remarks}

We will use the following version of the implicit function Theorem:

\begin{Theorem} [Implicit Function Theorem].
\label{Implicit Fct Thm}
Let $\pi: \CC^2 \times \calA \rr \calA$ be a real-analytic map defined on a 
neighborhood of $(x_0,\phi_0)\in \RR^2\times\calA_\RRs$. 
We let $T_0=D_\phi \pi(x_0,\phi_0)$
 denote the derivative of this map with respect to $\phi$.
Suppose that
$\phi_0=\pi(x_0,\phi_0)\in\calA_\RRs$ and that the spectral radius 
of the derivative,
$\rho(T_0)$, 
 is strictly smaller than $1$. Then
there exists  a neighborhood $U\subset \CC^2$
 of $x_0$ and a real-analytic map (unique if $U$ is small enough),
$x\in U \mapsto \phi(x)\in \calA$, for which $\phi_0=\phi(x_0)$,
$\phi(x)=\pi(x,\phi(x))$ 
and $\rho(D_\phi\pi(x,\phi(x)))<1$
for all $x\in U$.
\end{Theorem}
Proof: The map,
\[
    \Gamma(x,\phi) = (1-T_0)^{-1}
       (\pi(x,\phi)-\phi_0-T_0(\phi-\phi_0))+\phi_0,
\]
is real-analytic and verifies $\Gamma(x_0,\phi_0)=\phi_0$ and
$D_\phi\Gamma(x_0,\phi_0)=0$. We may therefore  find a neighborhood
$U$ of $x_0$ and a closed neighborhood $W$ of $\phi_0$ such that
$\Gamma$ is a uniform contraction on the real-analytic
sections, $U\rr W$. The fixed point $\phi(x)=\Gamma(x,\phi(x))$,
$x\in U$ is then itself a real-analytic section and has the desired
properties.  \Halmos\\

\begin{Lemma} The pressure function $P(s,\Lambda(\FF^t))$,
 extends to a real-analytic  function
$\PP(s,t)$, on an  open neighborhood $U^{s_0}\subset \CC^2$ of
$(s_0,0)$. 
\end{Lemma}
Proof:
By the above implicit function Theorem
there is a real-analytic map
 \[
    (s,t)\in U^{s_0} \mapsto h^*_{s,t} \in B(h^*_{s_0,0},\rho) \subset \calA.
 \]
defined in a neighborhood $U^{s_0}\subset W^{s_0}$ of $(s_0,0)$.
On this neighborhood we define as before,
$p_{s,t,\omega}=\ell(L_{s,F_{t,\omega}}h^*_{s,t}) \in \CC$.
 For fixed $\omega$ this function
is clearly analytic in $(s,t)\in U^{s_0}$.
Lemma
\ref{lemma neighborhood} applied to our fixed point shows that
when $(s,t)\in W^{s_0}$ and $\lambda_\omega=\|D\bff_{0,\omega}\|$ then
 \[
    \left|
      \frac{\lambda_\omega^s p_{s,t,\omega}}
      {\lambda_\omega^{s_0} p_{{s_0},0,\omega}} - 1 \right| \leq \frac23 .\]
This in turn implies,
      $|\log ({\lambda_\omega^s p_{s,t,\omega}})-
      \log (\lambda_\omega^{s_0} p_{{s_0},0,\omega})| \leq \log 3$.
Then also,
  \[ | \log p_{s,t,\omega} | \leq
       (|s|+s_0) \log \lambda_\omega + \log p_{s_0,0,\omega} + \log 3.\]
The right hand side is $\mu$-integrable 
(its integral is bounded by 
$(|s|+s_0) \EEE(\log \|D\bff_{0,\omega}\|) + P(s_0,\Lambda(\FF^0)) + \log 3$)
and therefore,
\[   \PP(s,t)=\int p_{s,t,\omega} \mu(d\omega), \ \ \ (s,t)\in U^{s_0} \]
is well-defined and yields
a real-analytic extension of the pressure.\Halmos\\

\begin{Theorem}
\label{part I}
Let $\tau$ be an ergodic transformation on $(\Omega,\mu)$.
Let $\FF^t=(\bff_{t,\omega})_{\omega\in\Omega}\in \EE_\Omega(K,U)$
be an analytic family verifying a  uniform $L$-Lipschitz condition and
with uniform bounded condition numbers, i.e.\  Assumption \ref{Assump} above.
 Then, almost surely, the
Hausdorff dimension of the random Julia set,
 $J(\FF^t_\omega)$, 
(\ref{random Julia set})
is independent of $\omega$ 
 and depends real-analytically on $t$.
\end{Theorem}

Proof:
Let $t\in \DD\cap \RR$.
We already know from Theorem 
\ref{Thm erg analytic}
that a.s.,
 $d(t)=\dimH \Lambda (\FF^t_\omega)$ is independent of $\omega$ and that
$\PP(d(t),t)=0$ whenever $(d(t),t)\in U^{s_0}$, $t\in \RR$.
By the previous Lemma, $\PP$ has a real-analytic extension and since
 $\frac{\partial \PP}{\partial s} (d(t),t) \leq \log \beta < 0$
for real $t$-values,
we may apply another implicit function theorem to $\PP$ and conclude
that there is an open neighborhood $V_0 \in \CC$ of $0$ and 
a real-analytic function $t\in V_0\mapsto (\hatd(t),t)\in U^{s_0}$
such that $\PP(\hatd(t),t)=0$ for all $t\in V_0$.
The function $\hatd(t)$ yields
 the desired real-analytic extension of
the dimension.
 \Halmos \\

\subsection{Parameter dependency of the measure}

Consider $\calM\equiv \calM(\Omega)$,
 the Banach space of complex measures on $\Omega$
in the variation norm. 
The set of probability measures,
$\calP\equiv \{\mu\in\calM: \mu\geq 0, \mu(\Omega)=1\}$,
forms a real affine subspace of $\calM$. Let
$\OO\subset \RR^n$ be an open subset.
\begin{Definition}
  We say that a family of probability measures, $p_\lambda$,
   $\lambda\in \OO$ is 
  real-analytic if there is a complex neighborhood 
  $\OO_\CCs\subset \CC^n$ of $\OO$, 
  such that
       \[
           \lambda\in \OO_\CCs \mapsto p_\lambda\in \calM
        \]
  is analytic.
\end{Definition}

\begin{Example}
  The Poission law,
  $p_\lambda(k)=e^{-\lambda} \frac{\lambda^k}{k!}$, $k\in \NN$
  is real-analytic in $\lambda\geq 0$. It has a complex
  extension to every $\lambda\in \CC$ with a variation norm
    \[
        \|p_\lambda\| = e^{|\lambda| - {\re} \lambda}
   \]
\end{Example}

\section {Proof of
Theorem \ref{Main Theorem}}
Let $(\Omega,\mu)=(\prod_\NNs \Upsilon,
\otimes_\NNs \nu)$ denote the (extension of the)
 direct product of probability spaces
and let $\tau$ be the shift on this space, i.e.\
$\tau(\bomega)=
(\bomega_2,\bomega_3,\ldots)$
for $\bomega=(\bomega_1,\bomega_2,\ldots)$.
With $f_{t,\omega}$ as in our Main Theorem we define
$\bff_{t,\bomega}=f_{t,\bomega_1}$
as  the random sequence of conformal maps.
We suppose that each individual measure $\nu_\lambda$ depends analytically
on a complex parameter $\lambda\in\DD$ (setting $\lambda=t$ we
obtain the statement in the Theorem).
The family $\FF^t=(\bff_{t,\bomega})_{\bomega\in\Omega}\in \EE_\Omega(K,U)$
verifies the conditions for Theorem
\ref{part I} and applying this for a (real) probability measure,
$\mu_\lambda$, $\lambda$ real,
yields part I and II of our Main Theorem,
except for the real-analyticity with respect
 to  the measure.
Let $\PP(s,t,\lambda)$ (for $\lambda$ real) denote the 
pressure obtained in that Theorem.

Going back to the Implicit Function Theorem, Theorem \ref{Implicit Fct Thm},
we may find a neighborhood $U^{s_0}\in \CC^2$ 
 of $(s_0,0)$ such that $D_h \pi_{s,t}(h^*_{s,t})$ has spectral radius
strictly smaller than one for $(s,t)\in U^{s_0}$. 
 Possibly shrinking
the neighborhood we may also  
find constants $C=C(s,t)<+\infty$, $\eta=\eta(s,t)<1$ 
and $0<\rho_1\leq \rho$ such that
the map $\pi_{s,t}^{(n)}: B(h^*,\rho_1) \rr B(h^*,\rho)$
is well-defined 
for all $n\geq 1$
and is a $C \eta^n$-Lipschitz contraction.

For $\eta<1$ we define,
$ D_\eta = \{\lambda \in \CC^n: \|p_\lambda\| < 1/\eta \}$, 
and then 
  \[ \D_0 = \{(s,t,\lambda): (s,t)\in U^{s_0},
       \lambda \in D_{\eta(s,t)} \}.\]
Given $(s,t)\in U^{s_0}$,
set $h^{(0)}=\bfone \equiv \pi^{(0)} \bfone \in B(h^*,\rho_1)$
 and then recursively,
 $h^{(n)} = \pi^{(n)} \bfone - \pi^{(n-1)} \bfone$,
$n\geq 1$. These differences have
norm smaller than $2 C\eta^n$.
Also $h^{(n)}_{\omega} = h^{(n)}_{\omega,\ldots,\tau^{n}\omega}$
 depends only on the
first $n$ iterates of $\omega$.
Integrating with respect to the analytic continuation
of our probability measure we see that
  \begin{equation}
   \sum_k |\ell(L_\omega h^{(k)}_{\omega_1 ...\omega_k}) 
                dp_\lambda(\omega_1)\ldots dp_\lambda(\omega_k)|
    \leq {\rm const} \sum_k (\|p_\lambda\| \eta)^k
  \label{real an ext}
   \end{equation}
which is finite when $\lambda\in D_\eta$. 
For  $\lambda$ real,
   \[
   \PP(s,t,\lambda)\equiv
        \int \ell(L_\omega h_{\tau \omega}) d\mu_\lambda(\omega)
   = \sum_k \ell(L_\omega h^{(k)}_{\omega_1 ...\omega_k}) 
                dp_\lambda(\omega_1)\ldots dp_\lambda(\omega_k)
 \]
and (\ref{real an ext}) shows that the right hand side
extends real-analytically on the domain
$(s,t,\lambda) \in \D_0$.
Using transversality of this extended  pressure function and once again
an Implicit Function Theorem we obtain  
Theorem \ref{Main Theorem}, part II,
including the real-analyticity with respect
to the measure.\Halmos\\

\begin{Remarks}
An alternative generalisation
 would be to pick the maps
$f_{t,\bomega}$ according to a Gibbs measure on a shift space over a finite
alphabet.   The Hausdorff dimension in this case depends real-analytically 
(and for the same reasons) 
upon the H\"older potential 
 defining the Gibbs state. 
This result does not, however,
cover our main example in the introduction.
\end{Remarks}

\section {Proof of Example \ref{example main}}
We define for $0\leq\rho<1$,  the complex annulus
$A_\rho=\{z\in\CC: \rho<|z|<1/\rho\}$ ($=\CC^*$ for $\rho=0$).
The conditions on parameters may be written as 
 \[  |a| + r \leq  \frac{k^2}{4}, \]
where $k$ is a constant $0< k <1$. We set 
$U=A_{k^2/2}$ and $K=\overline{A_{k/2}}$ which is a compact subset of $U$.

The maps under consideration, $f=z^{N+2}+c$, then belongs to 
$\EE(K,U)$.
The neighborhood $K_\Delta$ may be written as $\overline{A_{\kappa}}$ for
some $\kappa\in ]k^2/2,k/2[$.
 Conformal derivatives and usual derivatives are (smoothly) comparable
on $f^{-1}K_\Delta$ so we are allowed to replace conformal derivatives
by the standard Euclidean ones in the following.
 For $w=f(z)\in K_\Delta$ we have,
 \[  f'(z)= (N+2)z^{N+1}= (N+2) \frac{w-c}{z} ,\]
which is comparable to $N$ (because both $w$ and $z$ belongs to $K_\Delta$).
 Whence, the {\em b.a.l.d.}\ condition,
$\EEE(\log \|Df\|)<+\infty$, is equivalent to $\EEE(N)<+\infty$
which is clearly verified for a Poisson distribution of $N$.
Also, the condition numbers $\|Df\| \, \|1/Df\|$ are uniformly bounded
(this is in fact true for all maps $f\in \EE(K,U)$ for which 
$f^{-1}U$ is connected).
If we write
$f(z)=z^{N+2} + a + r \xi$, where $\xi$ is a  random variable
uniformly distributed in $\DD$ then we obtain an explicit
(real-) analytic parametrization of $f$ in terms of $a$ and $r$.

To see that a local inverse depends
uniformly Lipshitz in parameters consider e.g.~:
\[
 \frac
    {\partial f^{-1}}
    {\partial a} =
 -\frac
    {\partial f}
    {\partial a} /
    \frac{\partial f}
    {\partial z} =
    -\frac{1}{N+2}\; \frac{z}{w-c},
\]
which is uniformly bounded on $K_\Delta$.
Similarly,
\[
 \frac
    {\partial}
    {\partial a} \, \log f'\circ f^{-1}=
    \frac{N+1}{z}\; (- \frac{w-c}{(N+2)z})=
    \frac{N+1}{N+2}\; \frac{c-w}{z^2},
\]
which is again uniformly bounded, independent of the value of $N$
(but only just so !).
We are in the position to apply our Main Theorem and proving
the claims in  Example
\ref{example main}.\Halmos\\

\def\theequation{\Alph{section}.\arabic{equation}}
\appendix

\section{Removing the mixing condition}
\label{app removing}
Our mixing condition (C4) was convenient but not strictly necessary.
For completeness 
we will show  how to get rid of this condition.  
Our first reduction is to replace (C4) by topological
transitivity. This amounts to saying that there is $n_0=n_0(\delta)$
such that
      \begin{enumerate}
      \item[(C4')] $\displaystyle
              \bigcup_{k=0}^{n_0} f^{k}(B(x,\delta)\cap \Lambda)
                = \Lambda$.
      \end{enumerate}
Repeating the previous steps we see that 
(\ref{mM bounds}) is replaced by the inequality
 \[
     \max_{0\leq k\leq n_0} m_{n+k} 
            \geq (\lambda_1^{n_0} c_n)^{-s} M_n/2
 \]
from which the operator distortion bounds follow.
The proof of the lower bounds for the Hausdorff dimension
does not change and
in the upper bounds for the Box dimension the left hand side
of the inequality  (\ref{box lower})  is replaced by
$\sum_{0\leq j\leq n_0} L^{j+n_0} \chi_{B(x,r)}$ which leads
to the bound
  \[
              \sum_{i=1}^m (\diam\; B(x_i,2r))^s \leq
              {4^s}\gamma_2(s)
\sum_{0\leq j\leq n_0}\|L_s^{j}\|.
  \]

In the general situation we will replace $\Lambda$ by a
subset $\Lambda'$ 
of the same dimensions, which is $f$-invariant 
 and topologically transitive.
First, define a local pressure at $x\in \Lambda$ within $\Lambda$:
  \[
     \overline{P}_x(s,\Lambda) = \limsup_n \frac1n \log L_s^n \bfone_\Lambda(x).
  \]
\mbox{} From the very definitions it is clear that
$\barP_x(s,\Lambda)\leq \barP_{fx}(s,\Lambda)$.
Also if $x\in\Lambda$ and $y\in\Lambda$ are at a distance less than $\delta$
the ratio of $L^n_s \bfone (x)$ and $L^n_s \bfone(y)$ are 
sub-exponentially bounded in $n$. The local pressures at $x$ and $y$
are thus the same. Say that two points $x,y\in\Lambda$ are 
$\delta$-connected iff there is a finite sequence of points
  \[
    x_0=x,x_1,\ldots,x_n,x_{n+1}=y \subset \Lambda
  \]
for which $d(x_i,x_{i+1})<\delta$ for all $0\leq i\leq n$.
This partitions $\Lambda$ into
$\delta$-connected components $\Lambda=\Lambda_1\cup \ldots\cup \Lambda_m$.
Each $\Lambda_i$ is $\delta$-separated from its complement,
whence open and compact within $\Lambda$. Thus
there is a uniform bound on the number $N_\delta$ of intermediate
points needed to connect any $x$ and $y$ within the same component.

The partition is not Markovian. For example, for a connected hyperbolic 
Julia set there is only one $\delta$-connected component.
It does, however, enjoy some Markov like properties:
If $f \Lambda_i \cap \Lambda_j \neq \emptyset$ then
$f\Lambda_i \supset \Lambda_j$. To see this note that if
$x\in\Lambda_i$, $y=f(x_i)\in\Lambda_j$ and 
$v\in B(y,\delta)\subset \Lambda_j$ then there is (a unique)
$u\in B(x,\delta)\subset \Lambda_i$ for which $f(u)=v$
and thus $v\in \Lambda_i$.
We may introduce a transition matrix, $t_{ji}=1$ when $
f\Lambda_i\supset\Lambda_j$ and zero otherwise.
A partial ordering among the partition elements $\Lambda_i$ is then given
by 
  \[
    \Lambda_i \prec \Lambda_j \ \mbox{\ iff \ } \ 
    \exists n=n(i,j): t^n_{ji}\geq 1
  \]
and an equivalence relation 
  \[
    \Lambda_i \sim \Lambda_j \ \mbox{\ iff \ } \ 
    \Lambda_i \prec \Lambda_j  \ \mbox{\ and\ } \ 
    \Lambda_j \prec \Lambda_i .
  \]
The equivalence classes  provides a new partition of $\Lambda$:
  \[
     \Lambda = {\cal C}_1 \cup \ldots \cup {\cal C}_k 
  \]
which inherits the partial ordering from before.
Each equivalence class is topologically transitive and the
local pressures are constant on each class. Writing 
$P_i$ for the pressure on class $i$ we have 
$P_i \leq P_j$ for $i \prec j$.

Consider now the critical s-value
$s_\crit$ and let ${\cal C}_{i_0}$ be a class which is minimal
for the inherited partial ordering and such that the
local pressure vanishes for every point in this class
 $\overline{P}_x(s_\crit,\Lambda)=0$, $x\in{\cal C}_{i_0}$.
We denote by
  \[
     \Lambda'=\cap_{j\geq 0} f^{-j} {\cal C}_{i_0}
  \]
the corresponding 
 $f$ invariant subset of the class. This subset is
topologically transitive (clear) and
we claim that this set has Hausdorff and box dimensions that
agree and equal $s_\crit$.
 For this it  suffices to show that the pressure of that subset
$P(s_\crit,\Lambda',f)$ vanishes.\\

Write for $1\leq i \leq k$
\[
   N_i\phi = \chi_{{\cal C}_i} L_s \phi = 
L_s (\chi_{{\cal C}_i}\circ f \phi) .
\]
If ${\cal C}_i \prec {\cal C}_j$ and they are not equal then 
$N_i N_j\equiv 0$. Similarly, if 
 ${\cal C}_i$ and ${\cal C}_j$ are not related the
$N_i N_j=N_j N_i\equiv 0$. In either case we have
$(N_i+N_j)^n=N_i^n+ N_j N_i^{n-1}+ \cdots N_j^{n-1}N_i+N_j^n$ which
implies that the spectral radius of $N_i+N_j$ is the same as the
spectral radius of $N_j$.
Writing $L_s = \sum_i N_i$ it follows that the spectral
radius of $L_{s_\crit}$ must be the same as that of $N_{i_0}$.
 But this implies
precisely that $P(s_\crit,\Lambda',f)=0$.

\begin{Remarks}
\label{non finite measure}
We note that in this setting, even when distortions remain uniformly
bounded the Hausdorff measure need not be finite
 (essentially because the powers of a matrix of
 spectral radius one need not be 
bounded when the eigenvalue one is not simple).
\end{Remarks}

\end{document}